\documentclass[11pt]{amsart}
\usepackage{amsmath,amssymb,amsthm,mathrsfs,enumitem,url,tikz-cd}
\usepackage[hmargin=2cm,vmargin=3cm]{geometry}

\setlist[enumerate]{leftmargin=*}
\setlist[itemize]{labelindent=\parindent, leftmargin=*}

\numberwithin{equation}{section}

\theoremstyle{plain}
\newtheorem{thm}{Theorem}[section]
\newtheorem{lem}[thm]{Lemma}
\newtheorem{prop}[thm]{Proposition}
\newtheorem{cor}[thm]{Corollary}
\theoremstyle{definition}

\theoremstyle{remark}
\newtheorem{rem}[thm]{Remark}

\newcommand\Hom{\operatorname{Hom}}
\newcommand\Ind{\operatorname{Ind}}
\newcommand\Irr{\operatorname{Irr}}
\newcommand\JH{\operatorname{JH}}
\renewcommand\Re{\operatorname{Re}}
\newcommand\supp{\operatorname{supp}}

\newcommand\abs{\mathrm{abs}}
\newcommand\Ad{\mathrm{Ad}}
\newcommand\aut{\mathrm{aut}}
\newcommand\cusp{\mathrm{cusp}}
\newcommand\disc{\mathrm{disc}}
\newcommand\GL{\mathrm{GL}}
\newcommand\Mp{\mathrm{Mp}}
\newcommand\OO{\mathrm{O}}
\newcommand\PGL{\mathrm{PGL}}
\newcommand\SL{\mathrm{SL}}
\newcommand\SO{\mathrm{SO}}
\newcommand\Sp{\mathrm{Sp}}
\newcommand\Std{\mathrm{Std}}
\newcommand\Sym{\mathrm{Sym}}

\newcommand\1{\mathbf{1}}

\newcommand\A{\mathbb{A}}
\newcommand\C{\mathbb{C}}
\newcommand\F{\mathbb{F}}
\newcommand\Q{\mathbb{Q}}
\newcommand\R{\mathbb{R}}
\newcommand\V{\mathbb{V}}
\newcommand\Z{\mathbb{Z}}

\renewcommand\AA{\mathcal{A}}
\newcommand\TT{\mathcal{T}}

\newcommand\longhookrightarrow{\lhook\joinrel\longrightarrow}
\newcommand\longtwoheadrightarrow{\relbar\joinrel\twoheadrightarrow}
\newcommand\longhookleftarrow{\longleftarrow\joinrel\rhook}

\newcommand{\BIGOP}[1]{\mathop{\mathchoice%
{\raise-0.22em\hbox{\huge $#1$}}%
{\raise-0.05em\hbox{\Large $#1$}}{\hbox{\large $#1$}}{#1}}}
\newcommand{\BIGboxplus}{\mathop{\mathchoice%
{\raise-0.35em\hbox{\huge $\boxplus$}}%
{\raise-0.15em\hbox{\Large $\boxplus$}}{\hbox{\large $\boxtimes$}}{\boxtimes}}}

\title{The Shimura--Waldspurger correspondence for $\Mp_{2n}$}
\author{Wee Teck Gan}
\address{Department of Mathematics, National University of Singapore, 10 Lower Kent Ridge Road, Singapore 119076}
\email{matgwt@nus.edu.sg}
\author{Atsushi Ichino}
\address{Department of Mathematics, Kyoto University, Kitashirakawa Oiwake-cho, Sakyo-ku, Kyoto 606-8502, Japan}
\email{ichino@math.kyoto-u.ac.jp}

\begin{document}

\maketitle

\begin{abstract}
We generalize the Shimura--Waldspurger correspondence, which describes the generic part of the automorphic discrete spectrum of the metaplectic group $\Mp_2$, to the metaplectic group $\Mp_{2n}$ of higher rank.
To establish this, we transport Arthur's endoscopic classification of representations of the odd special orthogonal group $\SO_{2r+1}$ with $r \gg 2n$ by using a result of J.~S.~Li on global theta lifts in the stable range.
\end{abstract}

\section{\textbf{Introduction}}

In a seminal 1973 paper \cite{shimura}, Shimura revolutionized the study of half integral weight modular forms by establishing a lifting
\[
\begin{tikzcd}
 \{ \text{Hecke eigenforms of weight $k + \tfrac{1}{2}$ and level $\Gamma_0(4)$} \} \arrow[d] \\
 \{ \text{Hecke eigenforms of weight $2k$ and level $\SL_2(\Z)$} \}
\end{tikzcd}.
\]
He proved this by using Weil's converse theorem and a Rankin--Selberg integral for the standard $L$-function of a half integral weight modular form.
Subsequently, Niwa \cite{niwa} and Shintani \cite{shintani} explicitly constructed the Shimura lifting and its inverse by using theta series lifting.
From a representation theoretic viewpoint, Howe \cite{howe0} further developed the theory of theta series lifting by introducing the notion of reductive dual pairs, so that the lifting gives a correspondence between representations of one member of a reductive dual pair and those of the other member.
Then, in two influential papers \cite{w1,w2}, Waldspurger studied this correspondence in the framework of automorphic representations of the metaplectic group $\Mp_2$, which is a nonlinear two-fold cover of $\SL_2 = \Sp_2$.
Namely, he described the automorphic discrete spectrum of $\Mp_2$ precisely in terms of that of $\PGL_2 = \SO_3$ via the global theta lifts between $\Mp_2$ and (inner forms of) $\SO_3$.
Subordinate to this global result is the local Shimura correspondence, which is a classification of irreducible genuine representations of $\Mp_2$ in terms of that of $\SO_3$ and was also established by Waldspurger.
For an expository account of Waldspurger's result, the reader can consult a lovely paper of Piatetski-Shapiro \cite{ps} and a more recent paper \cite{g3} which takes advantage of 30 years of hindsight and machinery.

The goal of the present paper is to establish a similar description of the automorphic discrete spectrum of $\Mp_{2n}$, which is a nonlinear two-fold cover of $\Sp_{2n}$, in terms of that of $\SO_{2n+1}$.
As in the case of $\Mp_2$, it is natural to attempt to use the global theta lifts between $\Mp_{2n}$ and (inner forms of) $\SO_{2n+1}$ to relate these automorphic discrete spectra.
However, we encounter a difficulty.
For any irreducible cuspidal automorphic representation $\pi$ of $\Mp_{2n}$, there is an obstruction to the nonvanishing of its global theta lift to $\SO_{2n+1}$ given by the vanishing of the central $L$-value $L(\frac{1}{2}, \pi)$.
Thus, if we would follow Waldspurger's approach, then we would need the nonvanishing of the central $L$-value $L(\frac{1}{2}, \pi, \chi)$ twisted by some quadratic Hecke character $\chi$.
In the case of $\Mp_2$ (or equivalently $\PGL_2$), Waldspurger \cite{w2} proved the existence of such $\chi$ by exploiting Flicker's result \cite{flicker} on the correspondence between automorphic representations of $\GL_2$ and those of its two-fold cover.
Moreover, Jacquet \cite{jacquet} outlined a new proof of Waldspurger's nonvanishing result based on relative trace formulas, and Friedberg--Hoffstein \cite{fh} gave yet another proof and its extension to the case of $\GL_2$.
However, in the higher rank case, this seems to be a very difficult question in analytic number theory.
The main novelty of this paper is to overcome this inherent difficulty when $L(\frac{1}{2}, \pi) = 0$.

\subsection{Near equivalence classes}
\label{SS:near-eq-intro}

We now describe our results in more detail.
Let $F$ be a number field and $\A$ the ad\`ele ring of $F$.
Let $n$ be a positive integer.
We denote by $\Sp_{2n}$ the symplectic group of rank $n$ over $F$ and by $\Mp_{2n}(\A)$ the metaplectic two-fold cover of $\Sp_{2n}(\A)$:
\[
 1 \longrightarrow \{ \pm 1 \} \longrightarrow \Mp_{2n}(\A) \longrightarrow \Sp_{2n}(\A) \longrightarrow 1.
\]
Let $L^2(\Mp_{2n})$ be the subspace of $L^2(\Sp_{2n}(F) \backslash \Mp_{2n}(\A))$ on which $\{\pm 1\}$ acts as the nontrivial character, where we regard $\Sp_{2n}(F)$ as a subgroup of $\Mp_{2n}(\A)$ via the canonical splitting.
Then one of the basic problems is to understand the spectral decomposition of the unitary representation $L^2(\Mp_{2n})$ of $\Mp_{2n}(\A)$ and our goal is to establish a description of its discrete spectrum
\[
 L^2_{\disc}(\Mp_{2n})
\]
in the style of Arthur's conjecture formulated in \cite[Conjecture 25.1]{ggp}, \cite[\S 5.6]{g2}.

We first describe the decomposition of $L^2_{\disc}(\Mp_{2n})$ into near equivalence classes (which are coarser than equivalence classes) of representations.
Here we say that two irreducible representations $\pi = \bigotimes_v \pi_v$ and $\pi' = \bigotimes_v \pi_v'$ of $\Mp_{2n}(\A)$
are nearly equivalent if $\pi_v$ and $\pi'_v$ are equivalent for almost all places $v$ of $F$.
This decomposition will be expressed in terms of $A$-parameters defined as follows.
Consider a formal (unordered) finite direct sum
\begin{equation}
\label{eq:A-param}
 \phi = \bigoplus_i \phi_i \boxtimes S_{d_i}, 
\end{equation}
where
\begin{itemize}
\item $\phi_i$ is an irreducible self-dual cuspidal automorphic representation of $\GL_{n_i}(\A)$; 
\item $S_{d_i}$ is the unique $d_i$-dimensional irreducible representation of $\SL_2(\C)$.
\end{itemize}
We call $\phi$ an elliptic $A$-parameter for $\Mp_{2n}$ if
\begin{itemize}
\item $\sum_i n_i d_i  = 2n$;
\item if $d_i$ is odd, then $\phi_i$ is symplectic, i.e.~the exterior square $L$-function $L(s, \phi_i, \wedge^2)$ has a pole at $s=1$;
\item if $d_i$ is even, then $\phi_i$ is orthogonal, i.e.~the symmetric square $L$-function $L(s, \phi_i, \Sym^2)$ has a pole at $s=1$;
\item if $(\phi_i, d_i) = (\phi_j, d_j)$, then $i=j$.
\end{itemize}
If further $d_i = 1$ for all $i$, then we say that $\phi$ is generic.
We remark that the basic analytic properties of exterior and symmetric square $L$-functions were established in \cite{bf,js-aa,shahidi1,shahidi2} and \cite{bg,shahidi1,shahidi2}, respectively, and precisely one of $L(s, \phi_i, \wedge^2)$ and $L(s, \phi_i, \Sym^2)$ has a pole at $s=1$ since the Rankin--Selberg $L$-function
\[
 L(s, \phi_i \times \phi_i^\vee) = L(s, \phi_i \times \phi_i) = L(s, \phi_i, \wedge^2) L(s, \phi_i, \Sym^2)
\]
has a simple pole at $s=1$ by \cite[Proposition 3.6]{js} and both of the two $L$-functions do not vanish at $s=1$ by \cite[Theorem 5.1]{shahidi1}.
For each place $v$ of $F$,
let $\phi_v = \bigoplus_i \phi_{i,v} \boxtimes S_{d_i}$ be the localization of $\phi$ at $v$.
Here we regard $\phi_{i,v}$ as an $n_i$-dimensional representation of $L_{F_v}$ via the local Langlands correspondence \cite{langlands2,ht,henniart,scholze}, where
\[
 L_{F_v} =
 \begin{cases}
  \text{the Weil group of $F_v$} & \text{if $v$ is archimedean;} \\
  \text{the Weil--Deligne group of $F_v$} & \text{if $v$ is nonarchimedean.}
 \end{cases}
\]
Note that $\phi_v$ gives rise to an $A$-parameter $\phi_v : L_{F_v} \times \SL_2(\C) \rightarrow \Sp_{2n}(\C)$.
We associate to it an $L$-parameter $\varphi_{\phi_v} : L_{F_v} \rightarrow \Sp_{2n}(\C)$ by
\[
 \varphi_{\phi_v}(w) = \phi_v \!
 \left( w,
 \begin{pmatrix}
  |w|^{\frac{1}{2}} & \\
  & |w|^{-\frac{1}{2}}
 \end{pmatrix}
 \right).
\]

Our first result is:

\begin{thm}
\label{t:main1}
Fix a nontrivial additive character $\psi$ of $F \backslash \A$.
Then there exists a decomposition
\[
 L^2_{\disc}(\Mp_{2n}) = \bigoplus_{\phi} L^2_{\phi,\psi}(\Mp_{2n}),  
\]
where the direct sum runs over elliptic $A$-parameters $\phi$ for $\Mp_{2n}$ and $L^2_{\phi,\psi}(\Mp_{2n})$ is a full near equivalence class of irreducible representations $\pi$ in $L^2_{\disc}(\Mp_{2n})$ such that the $L$-parameter of $\pi_v$ (relative to $\psi_v$; see Remark \ref{r:unitarity} below) is $\varphi_{\phi_v}$ for almost all places $v$ of $F$.
\end{thm}

Thus, to achieve our goal, it remains to describe the decomposition of $L^2_{\phi,\psi}(\Mp_{2n})$ into equivalence classes of representations.
In this paper, we carry this out when $\phi$ is generic.

\begin{rem}
It immediately follows from Theorem \ref{t:main1} that $\Mp_{2n}$ has no embedded eigenvalues, i.e.~any family of eigenvalues of unramified Hecke algebras on the automorphic discrete spectrum of $\Mp_{2n}$ is distinct from that on the automorphic continuous spectrum of $\Mp_{2n}$.
This is an analog of Arthur's result \cite[Theorem 5]{a2011}, \cite{a}, which he needed to establish in the course of his proof of the classification of automorphic representations of orthogonal and symplectic groups.
However, in our case, we first establish the classification (with the help of theta lifts) and then deduce from it the absence of embedded eigenvalues.
\end{rem}

\begin{rem}
The dependence of $L^2_{\phi,\psi}(\Mp_{2n})$ on $\psi$ can be described as follows.
Let $\psi_a$ be another nontrivial additive character of $F \backslash \A$ given by $\psi_a(x) = \psi(ax)$ for some $a \in F^\times$.
If $\phi$ is of the form \eqref{eq:A-param}, we define another elliptic $A$-parameter $\phi_a$ by
\[
 \phi_a = \bigoplus_i (\phi_i \otimes (\chi_a \circ {\det})) \boxtimes S_{d_i},
\]
where $\chi_a$ is the quadratic automorphic character of $\A^\times$ associated to $F(\sqrt{a})/F$ by class field theory.
Then, for any irreducible summand $\pi$ of $L^2_{\phi,\psi_a}(\Mp_{2n})$, the $L$-parameter of $\pi_v$ (relative to $\psi_v$) is $\varphi_{\phi_{a,v}}$ for almost all $v$.
In particular, we have
\[
 L^2_{\phi,\psi_a}(\Mp_{2n}) = L^2_{\phi_a,\psi}(\Mp_{2n}).
\]
\end{rem}

\subsection{Local Shimura correspondence}

As in Waldspurger's result \cite{w1,w2}, our result will be expressed in terms of the local Shimura correspondence.
Fix a place $v$ of $F$ and assume for simplicity that $v$ is nonarchimedean.
For the moment, we omit the subscript $v$ from the notation, so that $F$ is a nonarchimedean local field of characteristic zero.
Let $\Irr \Mp_{2n}$ be the set of equivalence classes of irreducible genuine representations of the metaplectic group $\Mp_{2n}$ over $F$.
Then the local Shimura correspondence is a classification of $\Irr \Mp_{2n}$ in terms of $\Irr \SO(V)$, where $V$ is a $(2n+1)$-dimensional quadratic space over $F$.

To be precise, recall that there are precisely two such quadratic spaces with trivial discriminant (up to isometry).
We denote them by $V^+$ and $V^-$ so that $\SO(V^+)$ is split over $F$.
In \cite{gs}, the first-named author and Savin showed that for any nontrivial additive character $\psi$ of $F$, there exists a bijection 
\[
 \theta_\psi : \Irr \Mp_{2n} \longleftrightarrow \Irr \SO(V^+) \sqcup \Irr \SO(V^-)
\]
satisfying natural properties:
\begin{itemize}
 \item $\theta_\psi$ preserves the square-integrability and the temperedness of representations;
 \item $\theta_\psi$ is compatible with the theory of $R$-groups (modulo square-integrable representations);
 \item $\theta_\psi$ is compatible with the Langlands classification (modulo tempered representations);
 \item $\theta_\psi$ preserves the genericity of tempered representations;
 \item $\theta_\psi$ preserves various invariants such as $L$- and $\epsilon$-factors, Plancherel measures, and formal degrees \cite{gi1}.
\end{itemize}
Note that an analogous result in the archimedean case was proved by Adams--Barbasch \cite{ab1, ab2} more than 20 years ago.

On the other hand, the local Langlands correspondence, established by Arthur \cite{a} for $\SO(V^+)$ and by M{\oe}glin--Renard \cite{mr4} for $\SO(V^-)$, gives a partition
\[
 \Irr \SO(V^+) \sqcup \Irr \SO(V^-) = \bigsqcup_\phi \Pi_{\phi}(\SO(V^\pm)),
\]
where the disjoint union runs over equivalence classes of $L$-parameters $\phi:L_F \rightarrow \Sp_{2n}(\C)$ and $\Pi_{\phi}(\SO(V^\pm))$ is the associated Vogan $L$-packet equipped with a bijection
\[ 
 \Pi_{\phi}(\SO(V^\pm)) \longleftrightarrow \hat{S}_{\phi},
\]
where $S_{\phi}$ is the component group of the centralizer of the image of $\phi$ in $\Sp_{2n}(\C)$ and $\hat{S}_{\phi}$ is the group of characters of $S_{\phi}$.
Composing this with the local Shimura correspondence, we obtain a local Langlands correspondence for $\Mp_{2n}$:
\[
 \Irr \Mp_{2n} = \bigsqcup_{\phi} \Pi_{\phi,\psi}(\Mp_{2n})
\] 
with
\begin{equation}
\label{eq:bij-mp-intro}
 \Pi_{\phi,\psi}(\Mp_{2n}) \longleftrightarrow \hat{S}_\phi, 
\end{equation}
which depends on the choice of $\psi$, and which inherits various properties of the local Langlands correspondence for $\SO(V^\pm)$.
We remark that the dependence of the $L$-packet $\Pi_{\phi,\psi}(\Mp_{2n})$ and the bijection \eqref{eq:bij-mp-intro} on $\psi$ is described in \cite[Theorem 12.1]{gs}.
We also emphasize that the $L$-packet $\Pi_{\phi,\psi}(\Mp_{2n})$ satisfies \emph{endoscopic character relations} (see \S \ref{ss:endoscopy} below).

\subsection{Multiplicity formula for $\Mp_{2n}$}
\label{ss:mult-mp}

Suppose again that $F$ is a number field.
We now describe the multiplicity of any representation of $\Mp_{2n}(\A)$ in $L^2_{\phi,\psi}(\Mp_{2n})$ when $\phi$ is generic, i.e.~$\phi$ is a multiplicity-free sum
\[
 \phi = \bigoplus_i \phi_i
\]
of irreducible symplectic cuspidal automorphic representations $\phi_i$ of $\GL_{n_i}(\A)$.
We formally associate to $\phi$ a free $\Z/2\Z$-module
\[
 S_\phi = \bigoplus_i (\Z/2\Z) a_i
\]
with a basis $\{ a_i \}$, where $a_i$ corresponds to $\phi_i$.
We call $S_\phi$ the global component group of $\phi$.
For any place $v$ of $F$, this gives rise to a local $L$-parameter $\phi_v:L_{F_v} \rightarrow \Sp_{2n}(\C)$ together with a canonical map $S_\phi \rightarrow S_{\phi_v}$. 
Thus, we obtain a compact group $S_{\phi, \A} = \prod_v S_{\phi_v}$ equipped with the diagonal map $\Delta : S_\phi \rightarrow S_{\phi, \A}$.
Let $\hat{S}_{\phi, \A} = \bigoplus_v \hat{S}_{\phi_v}$ be the group of continuous characters of $S_{\phi, \A}$.
For any $\eta = \bigotimes_v \eta_v \in \hat{S}_{\phi, \A}$, we may form an irreducible genuine representation 
\[
 \pi_{\eta} = \bigotimes_v \pi_{\eta_v}
\]
of $\Mp_{2n}(\A)$, where $\pi_{\eta_v} \in \Pi_{\phi_v,\psi_v}(\Mp_{2n})$ is the representation of $\Mp_{2n}(F_v)$ associated to $\eta_v \in \hat{S}_{\phi_v}$.
Note that $\eta_v$ is trivial and $\pi_{\eta_v}$ is unramified for almost all $v$.
Also, by the Ramanujan conjecture for general linear groups, $\pi_\eta$ is expected to be tempered.
Finally, we define a quadratic character $\epsilon_\phi$ of $S_{\phi}$ by setting
\[
 \epsilon_\phi (a_i) = \epsilon(\tfrac{1}{2}, \phi_i),
\]
where $\epsilon(\frac{1}{2}, \phi_i) \in \{ \pm 1 \}$ is the root number of $\phi_i$.

Our second result (under the hypothesis that Arthur's result \cite{a} extends to the case of nonsplit odd special orthogonal groups; see \S \ref{ss:near-eq-so} and \S \ref{ss:mult-global} below for more details) is:

\begin{thm}
\label{t:main2}
Let $\phi$ be a generic elliptic $A$-parameter for $\Mp_{2n}$.
Then we have
\[
 L^2_{\phi,\psi}(\Mp_{2n}) \cong \bigoplus_{\eta \in \hat{S}_{\phi, \A}} m_\eta \pi_\eta,
\]
where
\[
 m_\eta = 
 \begin{cases}
  1 & \text{if $\Delta^* \eta = \epsilon_\phi$;} \\
  0 & \text{otherwise.}
 \end{cases}
\]
\end{thm}

\begin{rem}
The description of the automorphic discrete spectrum of $\Mp_{2n}$ is formally similar to that of $\SO_{2n+1}$, except that the condition $\Delta^* \eta = \epsilon_\phi$ is replaced by $\Delta^* \eta = \1$ in the case of $\SO_{2n+1}$.
\end{rem}

As an immediate consequence of Theorems \ref{t:main1} and \ref{t:main2}, we obtain the following generalization of Waldspurger's result \cite[p.~131]{w1}.

\begin{cor}
The generic part of $L^2_\disc(\Mp_{2n})$ (which is defined as $\bigoplus_\phi L^2_{\phi,\psi}(\Mp_{2n})$, where the direct sum runs over generic elliptic $A$-parameters $\phi$ for $\Mp_{2n}$) is multiplicity-free.
\end{cor}

\subsection{Idea of the proof}

The main ingredients in the proof of Theorems \ref{t:main1} and \ref{t:main2} are:
\begin{itemize}
\item 
Arthur's endoscopic classification \cite{a} (which relies on, among other things, the fundamental lemma proved by Ng\^o \cite{ngo} and the stabilization of the twisted trace formula established by M{\oe}glin--Waldspurger \cite{mw1,mw2});
\item
a result of J.~S.~Li \cite{li3} on global theta lifts in the stable range.
\end{itemize}
It is natural to attempt to transport Arthur's result for the odd special orthogonal group $\SO_{2n+1}$ to the metaplectic group $\Mp_{2n}$ by using global theta lifts.
However, as explained above, the difficulty arises when the central $L$-value vanishes.

To circumvent this difficulty, we consider the theta lift between $\Mp_{2n}$ and $\SO_{2r+1}$ with $r \gg 2n$, i.e.~the one in the stable range.
More precisely, let $\pi$ be an irreducible summand of $L^2_\disc(\Mp_{2n})$.
If $\pi$ is cuspidal, then by the Rallis inner product formula, the global theta lift $\theta_\psi(\pi)$ to $\SO_{2r+1}(\A)$ is always nonzero and square-integrable.
Even if $\pi$ is not necessarily cuspidal (so that the Rallis inner product formula is not available), J.~S.~Li \cite{li3} has developed a somewhat unconventional method for lifting $\pi$ to an irreducible summand $\theta_\psi(\pi)$ of $L^2_\disc(\SO_{2r+1})$, which is given by the spectral decomposition of the space of theta functions on $\SO_{2r+1}(\A)$ associated to the reductive dual pair $(\Mp_{2n},\SO_{2r+1})$.
(Note that such theta functions are square-integrable by the stable range condition.)
Then Arthur's result attaches an elliptic $A$-parameter $\phi'$ to $\theta_\psi(\pi)$, which turns out to be of the form 
\[
 \phi' = \phi \oplus S_{2r-2n}
\]
for some elliptic $A$-parameter $\phi$ for $\Mp_{2n}$ (see Proposition \ref{p:A-param}).
We now define the $A$-parameter of $\pi$ as $\phi$, which proves Theorem \ref{t:main1}.

To prove Theorem \ref{t:main2}, we apply Arthur's result to the near equivalence class $L^2_{\phi'}(\SO_{2r+1})$ and transport its local-global structure to $L^2_{\phi,\psi}(\Mp_{2n})$.
For this, we need the following multiplicity preservation: if
\[
 L^2_{\phi,\psi}(\Mp_{2n}) \cong \bigoplus_\pi m_\pi \pi,
\]
then
\[
 L^2_{\phi'}(\SO_{2r+1}) \cong \bigoplus_\pi m_\pi \theta_\psi(\pi).
\]
Since the above result of J.~S.~Li amounts to the theta lift from $\Mp_{2n}$ to $\SO_{2r+1}$, we need the theta lift from $\SO_{2r+1}$ to $\Mp_{2n}$ in the opposite direction.
In fact, J.~S.~Li \cite{li3} has also developed a method which allows us to lift an irreducible summand $\sigma$ of $L^2_{\phi'}(\SO_{2r+1})$ (which is no longer cuspidal so that the conventional method does not work) to an irreducible subrepresentation $\theta_\psi(\sigma)$ of the space of automorphic forms on $\Mp_{2n}(\A)$, which is realized as a Fourier--Jacobi coefficient of $\sigma$ (as in the automorphic descent; see \cite{grs99, grs02}).
However, we do not know a priori that $\theta_\psi(\sigma)$ occurs in $L^2_\disc(\Mp_{2n})$.
The key innovation in this paper is to show that $\theta_\psi(\sigma)$ is cuspidal if $\phi$ is generic (see Proposition \ref{p:tempered}).
From this, we can deduce the multiplicity preservation and hence obtain a multiplicity formula for $L^2_{\phi,\psi}(\Mp_{2n})$ when $\phi$ is generic (see Proposition \ref{p:mult-mp}).

However, there is still an issue: we do not know a priori that the local structure of $L^2_{\phi,\psi}(\Mp_{2n})$ transported from $L^2_{\phi'}(\SO_{2r+1})$ agrees with the one defined via the local Shimura correspondence.
In other words, we have to describe the local theta lift from $\SO_{2r+1}$ to $\Mp_{2n}$ in terms of the local Shimura correspondence (see Proposition \ref{p:key}).
This is the most difficult part in the proof of Theorem \ref{t:main2} and will be proved as follows.
\begin{itemize}
\item
We consider the theta lift of representations in the local $A$-packet $\Pi_{\phi'}(\SO_{2r+1})$ to $\Mp_{2n}$, where $\phi'$ is a local $A$-parameter of the form
\[
 \phi' = \phi \oplus S_{2r-2n} 
\]
for some local $L$-parameter $\phi$ for $\Mp_{2n}$.
For our global applications, we may assume that $\phi$ is almost tempered.
Then we can reduce the general case to the case of good $L$-parameters for smaller metaplectic groups, where we say that an $L$-parameter $\phi$ is good if any irreducible summand of $\phi$ is symplectic.
This will be achieved by using irreducibility of some induced representations (see Lemma \ref{l:irred}), which is due to M{\oe}glin \cite{m1,m2,m3,m4} in the nonarchimedean case and to M{\oe}glin--Renard \cite{mr} in the complex case, and which is proved in \cite{gi-mp-real} in the real case.
\item
If $\phi$ is good, then we appeal to a global argument.
As in our previous paper \cite{gi2}, we can find a global generic elliptic $A$-parameter $\varPhi$ such that $\varPhi_{v_0} = \phi$ for some $v_0$ and $\varPhi_v$ is non-good for all $v \ne v_0$, and then apply Arthur's multiplicity formula (viewed as a product formula) to extract information at $v_0$ from the knowledge at all $v \ne v_0$.
Strictly speaking, we need to impose more conditions on $\varPhi$ and the most crucial one is the nonvanishing of the central $L$-value $L(\frac{1}{2}, \varPhi)$ (see Corollary \ref{c:global}), which makes the argument more complicated than that in \cite{gi2}.
\end{itemize}

Finally, we remark that when $n=1$, our argument gives a new proof of the Shimura--Waldspurger correspondence for $\Mp_2$ which is independent of Waldspurger's result \cite[Th\'eor\`eme 4]{w2} on the nonvanishing of central $L$-values.

\subsection{Endoscopy for $\Mp_{2n}$}
\label{ss:endoscopy}

In \cite{wwli1,wwli2,wwli3}, W.-W.~Li has developed the theory of endoscopy for $\Mp_{2n}$ and has stabilized the elliptic part of the trace formula for $\Mp_{2n}$, which should yield a definition of local $L$-packets for $\Mp_{2n}$ satisfying endoscopic character relations.
In this paper, the local $L$-packets for $\Mp_{2n}$ are defined via the local Shimura correspondence and we do not know a priori that they satisfy the endoscopic character relations.
However, this was established by Adams \cite{adams} and Renard \cite{renard} in the real case.
Moreover, using a simple stable trace formula for $\Mp_{2n}$ established by W.-W.~Li and the main result of this paper as key inputs, C.~Luo \cite{luo}, a student of the first-named author, has recently proved the endoscopic character relations in the nonarchimedean case.

\subsection*{Acknowledgments}

We thank the referees for useful suggestions.
The first-named author is partially supported by an MOE Tier one grant R-146-000-228-114.
The second-named author is partially supported by JSPS KAKENHI Grant Numbers 26287003.

\subsection*{Notation}

If $F$ is a number field and $G$ is a reductive algebraic group defined over $F$, we denote by $\AA(G)$ the space of automorphic forms on $G(\A)$, where $\A$ is the ad\`ele ring of $F$.
If $G = \Mp_{2n}$, we understand that $\AA(G)$ consists only of genuine functions.
We denote by $\AA^2(G)$ and $\AA_\cusp(G)$ the subspaces of square-integrable automorphic forms and cusp forms on $G(\A)$, respectively.

If $F$ is a local field and $G$ is a reductive algebraic group defined over $F$, we denote by $\Irr G$ the set of equivalence classes of irreducible smooth representations of $G$, where we identify $G$ with its group of $F$-valued points $G(F)$.
If $G = \Mp_{2n}$, we understand that $\Irr G$ consists only of genuine representations.

For any irreducible representation $\pi$, we denote by $\pi^\vee$ its contragredient representation.
For any abelian locally compact group $S$, we denote by $\hat{S}$ the group of continuous characters of $S$.
For any positive integer $d$, we denote by $S_d$ the unique $d$-dimensional irreducible representation of $\SL_2(\C)$.

\section{\textbf{Some results of J.~S.~Li}}

In this section, we recall some results of J.~S.~Li \cite{li2, li3} on theta lifts and unitary representations of low rank which will play a crucial role in this paper.

\subsection{Metaplectic and orthogonal groups}

Let $F$ be either a number field or a local field of characteristic zero.
Let $n$ be a nonnegative integer and $\Sp_{2n}$ the symplectic group of rank $n$ over $F$.
(If $n=0$, we interpret $\Sp_0$ as the trivial group $\{ 1 \}$.)
If $F$ is local, we denote by $\Mp_{2n}$ the metaplectic two-fold cover of $\Sp_{2n}$:
\[
 1 \longrightarrow \{\pm 1\} \longrightarrow \Mp_{2n} \longrightarrow \Sp_{2n} \longrightarrow 1, 
\]
which is the unique (up to unique isomorphism) nonsplit two-fold topological central extension of $\Sp_{2n}$ when $n>0$ and $F \ne \C$, and which is the trivial extension (i.e.~$\Mp_{2n} \cong \Sp_{2n} \times \{ \pm 1 \}$) when $n=0$ or $F = \C$.
We may realize $\Mp_{2n}$ as the set $\Sp_{2n} \times \{ \pm 1 \}$ with multiplication law determined by Ranga Rao's two-cocycle \cite{rangarao}.
If $F$ is global, we denote by $\Mp_{2n}(\A)$ the metaplectic two-fold cover of $\Sp_{2n}(\A)$:
\[
 1 \longrightarrow \{ \pm 1\} \longrightarrow \Mp_{2n}(\A) \longrightarrow \Sp_{2n}(\A) \longrightarrow 1.
\]
This cover splits over $\Sp_{2n}(F)$ canonically.

Let $V$ be a quadratic space over $F$, i.e.~a finite-dimensional vector space over $F$ equipped with a nondegenerate quadratic form $q:V \rightarrow F$.
We assume that $V$ is odd-dimensional:
\[
 \dim V = 2r+1,
\]
and that the discriminant of $q$ is trivial.
We say that $V$ is split if $V$ is isometric to the orthogonal direct sum $\mathbb{H}^r \oplus F$, where $\mathbb{H}$ is the hyperbolic plane and $F$ is equipped with a quadratic form $q(x) = x^2$.
We denote by $\OO(V)$ the orthogonal group of $V$ and by $\SO(V) = \OO(V)^0$ the special orthogonal group of $V$.
Note that
\begin{equation}
\label{E:O(V)SO(V)}
 \OO(V) = \SO(V) \times \{ \pm 1\}. 
\end{equation}
We write
\[
 \SO(V) = \SO_{2r+1}
\]
if $V$ is split.
Then $\SO_{2r+1}$ is split over $F$.
If $F$ is local, we denote by $\varepsilon(V) \in \{ \pm 1\}$ the normalized Hasse--Witt invariant of $V$ (see \cite[pp.~80--81]{scharlau}).
In particular, $\varepsilon(V) = 1$ if $V$ is split.

\subsection{Theta lifts}

Suppose first that $F$ is local and fix a nontrivial additive character $\psi$ of $F$.
Recall that with the above notation, the pair $(\Mp_{2n}, \OO(V))$ is an example of a reductive dual pair.
Let $\omega_\psi$ be the Weil representation of $\Mp_{2n} \times \OO(V)$ with respect to $\psi$.
For any irreducible genuine representation $\pi$ of $\Mp_{2n}$, the maximal $\pi$-isotypic quotient of $\omega_\psi$ is of the form
\[
 \pi \boxtimes \Theta_\psi(\pi)
\]
for some representation $\Theta_\psi(\pi)$ of $\OO(V)$.
Then, by the Howe duality \cite{howe3, w3, gt2}, $\Theta_\psi(\pi)$ has a unique irreducible (if nonzero) quotient $\theta_\psi(\pi)$.
We regard $\theta_\psi(\pi)$ as a representation of $\SO(V)$ by restriction.
By \eqref{E:O(V)SO(V)}, $\theta_\psi(\pi)$ remains irreducible if nonzero.

Similarly, for any irreducible representation $\tilde\sigma$ of $\OO(V)$, 
we define a representation $\Theta_\psi(\tilde\sigma)$ of $\Mp_{2n}$ with its unique irreducible (if nonzero) quotient $\theta_\psi(\tilde\sigma)$.
We now assume that $r \ge n$.
Let $\sigma$ be an irreducible representation of $\SO(V)$.
By the conservation relation \cite{sz},
there exists at most one extension $\tilde\sigma$ of $\sigma$ to $\OO(V)$ such that $\theta_\psi(\tilde\sigma)$ is nonzero,
in which case we write $\theta_\psi(\tilde\sigma) = \theta_\psi(\sigma)$.
(We interpret $\theta_\psi(\sigma)$ as zero if such $\tilde{\sigma}$ does not exist.)

Suppose next that $F$ is global and fix a nontrivial additive character $\psi$ of $F \backslash \A$.
Let $\pi = \bigotimes_v \pi_v$ be an abstract irreducible genuine representation of $\Mp_{2n}(\A)$.
Assume that the theta lift $\theta_{\psi_v}(\pi_v)$ of $\pi_v$ to $\SO(V)(F_v)$ is nonzero for all places $v$ of $F$.
Then $\theta_{\psi_v}(\pi_v)$ is irreducible for all $v$ and is unramified for almost all $v$.
Hence we may define an abstract irreducible representation
\[
 \theta_\psi^{\abs}(\pi) = \bigotimes_v \theta_{\psi_v}(\pi_v)
\] 
of $\SO(V)(\A)$.
We call $\theta_\psi^{\abs}(\pi)$ the abstract theta lift of $\pi$ to $\SO(V)(\A)$.
On the other hand, if $\pi$ is an irreducible genuine cuspidal automorphic representation of $\Mp_{2n}(\A)$,
then we may define its global theta lift $\Theta_\psi^{\aut}(\pi)$ as the subspace of $\AA(\SO(V))$ spanned by all automorphic forms of the form
\[
 \theta(f,\varphi)(h) = \int_{\Sp_{2n}(F) \backslash \Mp_{2n}(\A)} \theta(f)(g,h) \, \overline{\varphi(g)} \, dg
\]
for $f \in \omega_\psi$ and $\varphi \in \pi$.
Here $\omega_\psi$ is the Weil representation of $\Mp_{2n}(\A) \times \OO(V)(\A)$ with respect to $\psi$ and $\theta(f)$ is the theta function associated to $f$, which is a slowly increasing function on $\Sp_{2n}(F) \backslash \Mp_{2n}(\A) \times \OO(V)(F) \backslash \OO(V)(\A)$.
If $\Theta_\psi^{\aut}(\pi)$ is nonzero and contained in $\AA^2(\SO(V))$, then $\Theta_\psi^{\aut}(\pi)$ is irreducible and 
\[ 
 \Theta_\psi^{\aut}(\pi) \cong \theta_\psi^{\abs}(\pi)
\]
by \cite[Corollary 7.1.3]{kr}.

\subsection{Unitary representations of low rank}
\label{ss:low-rank}

The notion of rank for unitary representations was first introduced by Howe \cite{howe2} in the case of symplectic groups and was extended by J.~S.~Li \cite{li2} to the case of classical groups.
Following \cite[\S 4]{li2}, we say that an irreducible unitary representation of $\SO_{2r+1}$ is of low rank if its rank (which is necessarily even) is less than $r-1$.
Such representations are obtained by theta lifts as follows.

Assume that $2n < r-1$.
In particular, the reductive dual pair $(\Mp_{2n},\SO_{2r+1})$ is in the stable range (see \cite[Definition 5.1]{li1}).
If $F$ is local, then for any irreducible genuine representation $\pi$ of $\Mp_{2n}$, its theta lift $\theta_\psi(\pi)$ to $\SO_{2r+1}$ is nonzero.
Moreover, if $\pi$ is unitary, then so is $\theta_\psi(\pi)$ by \cite{li1}.
In \cite{li2}, J.~S.~Li showed that this theta lift provides a bijection
\begin{equation}
\label{eq:li-bij}
\begin{tikzcd}
 \{ \text{irreducible genuine unitary representations of $\Mp_{2n}$} \}
 \times
 \{ \text{quadratic characters of $F^\times$} \} \arrow[d] \\
 \{ \text{irreducible unitary representations of $\SO_{2r+1}$ of rank $2n$} \} \arrow[u]
\end{tikzcd},
\end{equation}
which sends a pair $(\pi, \chi)$ in the first set to a representation $\theta_\psi(\pi) \otimes (\chi \circ \nu)$ of $\SO_{2r+1}$, where $\nu : \SO_{2r+1} \rightarrow F^\times/(F^\times)^2$ is the spinor norm.

This result has a global analog.
Let $F$ be a number field and $\sigma = \bigotimes_v \sigma_v$ an irreducible unitary representation of $\SO_{2r+1}(\A)$ which occurs as a subrepresentation of $\AA(\SO_{2r+1})$.
Then, by \cite{howe1}, \cite[Lemma 3.2]{li3}, the following are equivalent:
\begin{itemize}
 \item $\sigma$ is of rank $2n$;
 \item $\sigma_v$ is of rank $2n$ for all $v$;
 \item $\sigma_v$ is of rank $2n$ for some $v$.
\end{itemize}
Suppose that $\sigma$ satisfies the above equivalent conditions.
Then, for any $v$, there exist a unique irreducible genuine unitary representation $\pi_v$ of $\Mp_{2n}(F_v)$ and a unique quadratic character $\chi_v$ of $F_v^\times$ such that
\[
 \sigma_v \cong \theta_{\psi_v}(\pi_v) \otimes (\chi_v \circ \nu).
\]
By \cite[Proposition 5.7]{li3}, $\chi_v$ is unramified for almost all $v$ and the abstract character $\chi = \bigotimes_v \chi_v$ of $\A^\times$ is in fact automorphic.
This implies that $\theta_{\psi_v}(\pi_v)$ and hence $\pi_v$ are unramified for almost all $v$.
Hence we may define an abstract representation $\pi = \bigotimes_v \pi_v$ of $\Mp_{2n}(\A)$, so that
\[
 \sigma \cong \theta_\psi^{\abs}(\pi) \otimes (\chi \circ \nu).
\]

\subsection{Some inequalities}

Finally, we recall a result of J.~S.~Li \cite{li3} which allows us to lift square-integrable (but not necessarily cuspidal) automorphic representations of $\Mp_{2n}(\A)$ to $\SO_{2r+1}(\A)$.
For any irreducible genuine representation $\pi$ of $\Mp_{2n}(\A)$, we define its multiplicities $m(\pi)$ and $m_\disc(\pi)$ by
\begin{align*}
 m(\pi) & = \dim \Hom_{\Mp_{2n}(\A)}(\pi, \AA(\Mp_{2n})), \\
 m_\disc(\pi) & = \dim \Hom_{\Mp_{2n}(\A)}(\pi, \AA^2(\Mp_{2n})).
\end{align*}
Obviously, $m_\disc(\pi) \leq m(\pi)$.
Likewise, if $\sigma$ is an irreducible representation of $\SO_{2r+1}(\A)$, we have its multiplicities $m(\sigma)$ and $m_\disc(\sigma)$.

\begin{thm}[J.~S.~Li \cite{li3}]
\label{T:inequalities}
Assume that $2n < r-1$.
Let $\pi$ be an irreducible genuine unitary representation of $\Mp_{2n}(\A)$ and $\theta_\psi^{\abs}(\pi)$ its abstract theta lift to $\SO_{2r+1}(\A)$.
Then we have
\[
 m_\disc(\pi) \le m_\disc(\theta_\psi^{\abs}(\pi)) \le m(\theta_\psi^{\abs}(\pi)) \le m(\pi).
\]
\end{thm}

\section{\textbf{Near equivalence classes and $A$-parameters}}

In this section, we attach an $A$-parameter to each near equivalence class in $L^2_\disc(\Mp_{2n})$.

\subsection{The automorphic discrete spectrum of $\SO(V)$}
\label{ss:near-eq-so}

We first describe the automorphic discrete spectrum
\[
 L^2_{\disc}(\SO(V)) = L^2_{\disc}(\SO(V)(F) \backslash \SO(V)(\A)),
\]
where $V$ is a $(2r+1)$-dimensional quadratic space over a number field $F$ with trivial discriminant.
If $\SO(V)$ is split over $F$, then it follows from Arthur's result \cite{a} that
\[
 L^2_{\disc}(\SO(V)) = \bigoplus_{\phi} L^2_{\phi}(\SO(V)),
\]
where the direct sum runs over elliptic $A$-parameters $\phi$ for $\SO(V)$ (or equivalently those for $\Mp_{2r}$) and $L^2_{\phi}(\SO(V))$ is a full near equivalence class of irreducible representations $\sigma$ in $L^2_{\disc}(\SO(V))$ such that the $L$-parameter of $\sigma_v$ is $\varphi_{\phi_v}$ for almost all places $v$ of $F$.
Even if $\SO(V)$ is not necessarily split over $F$, this decomposition is expected to hold.

\subsection{Attaching $A$-parameters}

Let $C$ be a near equivalence class in $L^2_{\disc}(\Mp_{2n})$.
Then $C$ gives rise to a collection of $L$-parameters
\[
 \varphi_v : L_{F_v} \longrightarrow \Sp_{2n}(\C)
\]
for almost all $v$ such that for any irreducible summand $\pi$ of $C$, the $L$-parameter of $\pi_v$ (relative to $\psi_v$) is $\varphi_v$ for almost all $v$.
Here $\psi$ is the fixed nontrivial additive character of $F \backslash \A$.
Then we have:

\begin{prop}
\label{p:A-param}
There exists a unique elliptic $A$-parameter $\phi$ for $\Mp_{2n}$ such that $\varphi_{\phi_v} = \varphi_v$ for almost all $v$.
\end{prop}

\begin{proof}
To prove the existence of $\phi$, we fix an integer $r > 2n+1$ and consider the abstract theta lift from $\Mp_{2n}(\A)$ to $\SO_{2r+1}(\A)$.
Choose any irreducible summand $\pi$ of $C$.
Since $m_\disc(\pi) \ge 1$, we deduce from Theorem \ref{T:inequalities} that $m_\disc(\theta_\psi^\abs(\pi)) \ge 1$, i.e.~$\theta_\psi^{\abs}(\pi)$ occurs in $L^2_\disc(\SO_{2r+1})$.
Hence, as explained in \S \ref{ss:near-eq-so} above, Arthur's result \cite{a} attaches an elliptic $A$-parameter $\phi'$ to $\theta_\psi^{\abs}(\pi)$.

We show that $\phi'$ contains $S_{2r-2n}$ as a direct summand.
Consider the partial $L$-function $L^S(s, \theta_\psi^{\abs}(\pi))$ of $\theta_\psi^{\abs}(\pi)$ relative to the standard representation of $\Sp_{2r}(\C)$, where $S$ is a sufficiently large finite set of places of $F$.
If we write $\phi' = \bigoplus_i \phi_i \boxtimes S_{d_i}$ as in \eqref{eq:A-param}, then
\begin{equation}
\label{E:L1}  
 L^S(s, \theta_\psi^{\abs}(\pi)) = \prod_i \prod_{j=1}^{d_i} L^S \! \left( s + \frac{d_i+1}{2}-j, \phi_i \right).
\end{equation}
Note that $L^S(s, \phi_i)$ is holomorphic for $\Re s > 1$ for all $i$, and it has a pole at $s=1$ if and only if $\phi_i$ is the trivial representation of $\GL_1(\A)$.
On the other hand, by the local theta correspondence for unramified representations, the $L$-parameter of $\theta_{\psi_v}(\pi_v)$ is
\begin{equation}
\label{E:unram}
 \varphi_v \oplus \Big( \bigoplus_{j=1}^{2r-2n} |\cdot|^{r-n+\frac{1}{2}-j} \Big)
\end{equation}
for almost all $v$.
Hence 
\begin{equation}
\label{E:L2}
 L^S(s, \theta_\psi^{\abs}(\pi)) = L_\psi^S(s, \pi) \prod_{j=1}^{2r-2n} \zeta^S ( s+r-n+\tfrac{1}{2}- j),
\end{equation}
where $L^S_\psi(s, \pi)$ is the partial $L$-function of $\pi$ relative to $\psi$ and the standard representation of $\Sp_{2n}(\C)$.
Since $L^S_\psi(s, \pi)$ is holomorphic for $\Re s > n+1$ (see e.g.~\cite[Theorem 9.1]{y}), it follows from \eqref{E:L2} that $L^S(s, \theta_\psi^{\abs}(\pi))$ is holomorphic for $\Re s > r-n+\frac{1}{2}$ but has a pole at $s = r-n+\frac{1}{2}$.
This and \eqref{E:L1} imply that $\phi'$ contains $S_{2t}$ as a direct summand for some $t \ge r-n$.
If $t$ is the largest such integer, then $L^S(s, \theta_\psi^{\abs}(\pi))$ has a pole at $s = t+\frac{1}{2}$.
This forces $t=r-n$.

Thus, we may write
\[
 \phi' = \phi \oplus S_{2r-2n}
\]
for some elliptic $A$-parameter $\phi$.
This and \eqref{E:unram} imply that $\varphi_{\phi_v} = \varphi_v$ for almost all $v$.
Moreover, by the strong multiplicity one theorem \cite[Theorem 4.4]{js}, 
$\phi$ is uniquely determined by this condition.
\end{proof}

We now denote by $L^2_{\phi,\psi}(\Mp_{2n})$ the near equivalence class $C$, where $\phi$ is the $A$-parameter attached to $C$ by Proposition \ref{p:A-param}.
Then we have a decomposition
\[
 L^2_{\disc}(\Mp_{2n}) = \bigoplus_{\phi} L^2_{\phi,\psi}(\Mp_{2n}),
\]
where the direct sum runs over elliptic $A$-parameters $\phi$ for $\Mp_{2n}$.
Note that $L^2_{\phi,\psi}(\Mp_{2n})$ is possibly zero for some $\phi$.
This completes the proof of Theorem \ref{t:main1}.

\begin{rem}
\label{rem:embed}
Suppose that
\[
 L^2_{\phi,\psi}(\Mp_{2n}) \cong \bigoplus_{\pi} m_\pi \pi,
\]
where the direct sum runs over irreducible genuine representations $\pi$ of $\Mp_{2n}(\A)$ and $m_\pi$ is the multiplicity of $\pi$ in $L^2_{\phi,\psi}(\Mp_{2n})$.
If $r>2n+1$ and $\phi' = \phi \oplus S_{2r-2n}$, then it follows from Theorem \ref{T:inequalities} and the Howe duality that
\[
 L^2_{\phi'}(\SO_{2r+1}) \longhookleftarrow \bigoplus_{\pi} m_\pi \theta_\psi^{\abs}(\pi). 
\]
In Corollary \ref{c:mult-preserve} below, we show that the above embedding is an isomorphism if $\phi$ is generic.
\end{rem}

\begin{rem}
By Arthur's multiplicity formula \cite[Theorem 1.5.2]{a}, once we know that local $A$-packets are multiplicity-free, we have
\[
 m_\disc(\sigma) \le 1
\]
for all irreducible representations $\sigma$ of $\SO_{2r+1}(\A)$.
On the other hand, the multiplicity-freeness of local $A$-packets was proved by M{\oe}glin \cite{m1,m2,m3,m4} in the nonarchimedean case, and by M{\oe}glin \cite{m5} and M{\oe}glin--Renard \cite{mr} in the complex case, but is not fully known in the real case, though there has been progress by Arancibia--M{\oe}glin--Renard \cite{amr} and M{\oe}glin--Renard \cite{mr2,mr3}.
Hence, by Theorem \ref{T:inequalities}, $L^2_\disc(\Mp_{2n})$ is multiplicity-free at least when $F$ is totally imaginary.
\end{rem}

\section{\textbf{The case of generic elliptic $A$-parameters}}

In this section, we study the structure of $L^2_{\phi,\psi}(\Mp_{2n})$ for generic elliptic $A$-parameters $\phi$.

\subsection{A key equality}
\label{ss:key-eq}

For any irreducible genuine representation $\pi$ of $\Mp_{2n}(\A)$, we define the multiplicity $m_\cusp(\pi)$ by
\[
 m_\cusp(\pi) = \dim \Hom_{\Mp_{2n}(\A)}(\pi, \AA_\cusp(\Mp_{2n})).
\]
Obviously, $m_\cusp(\pi) \le m_\disc(\pi) \le m(\pi)$.

\begin{prop}
\label{p:tempered}
Let $\phi$ be a generic elliptic $A$-parameter for $\Mp_{2n}$.
Let $\pi$ be an irreducible genuine representation of $\Mp_{2n}(\A)$ such that the $L$-parameter of $\pi_v$ (relative to $\psi_v$) is $\phi_v$ for almost all $v$.
Then we have
\[
 m_\cusp(\pi) = m_\disc(\pi) = m(\pi).
\]
\end{prop}

\begin{proof}
It suffices to show that for any realization $\mathcal{V} \subset \AA(\Mp_{2n})$ of $\pi$, we have $\mathcal{V} \subset \AA_\cusp(\Mp_{2n})$.
Suppose on the contrary that $\mathcal{V} \not\subset \AA_\cusp(\Mp_{2n})$ for some such $\mathcal{V}$.
Then the image $\mathcal{V}_P$ of $\mathcal{V}$ under the constant term map
\[
 \AA(\Mp_{2n}) \longrightarrow \AA_P(\Mp_{2n})
\]
is nonzero for some proper parabolic subgroup $P$ of $\Sp_{2n}$.
Here $\AA_P(\Mp_{2n})$ is the space of genuine automorphic forms on $N(\A) M(F) \backslash \Mp_{2n}(\A)$, where $M$ and $N$ are a Levi component and the unipotent radical of $P$, respectively.
Assume that $P$ is minimal with this property, so that $\mathcal{V}_P$ is contained in the space of cusp forms in $\AA_P(\Mp_{2n})$.
Then, as explained in \cite[p.~205]{langlands1}, $\pi$ is a subrepresentation of
\[
 \Ind^{\Mp_{2n}(\A)}_{\tilde{P}(\A)}(\rho)
\]
for some irreducible cuspidal automorphic representation $\rho$ of $\tilde{M}(\A)$, where $\tilde{P}(\A)$ and $\tilde{M}(\A)$ are the preimages of $P(\A)$ and $M(\A)$ in $\Mp_{2n}(\A)$, respectively.
If $M \cong \prod_i \GL_{k_i} \times \Sp_{2n_0}$ with $\sum_i k_i + n_0 = n$, then $\rho$ is of the form
\[
 \rho \cong \Big( \bigotimes_i \tilde{\tau}_{i,\psi} \Big) \otimes \pi_0
\] 
for some irreducible cuspidal automorphic representations $\tau_i$ and $\pi_0$ of $\GL_{k_i}(\A)$ and $\Mp_{2n_0}(\A)$, respectively.
Here, as in \cite[\S 2.6]{gi1}, we define the twist $\tilde{\tau}_{i,\psi} = \tau_i \otimes \chi_\psi$ of $\tau_i$ by the genuine quartic automorphic character $\chi_\psi$ of the two-fold cover of $\GL_{k_i}(\A)$.
By Proposition \ref{p:A-param}, $\pi_0$ has a weak lift $\tau_0$ to $\GL_{2n_0}(\A)$.
Hence $\pi$ has a weak lift to $\GL_{2n}(\A)$ of the form
\begin{equation}
\label{E:weak-lift1}
 \Big( \BIGboxplus_i (\tau_i \boxplus \tau_i^{\vee}) \Big) \boxplus \tau_0,
\end{equation}
where $\boxplus$ denotes the isobaric sum.

On the other hand, if we write $\phi = \bigoplus_i \phi_i$ for some pairwise distinct irreducible symplectic cuspidal automorphic representations $\phi_i$ of $\GL_{n_i}(\A)$, then $\pi$ has a weak lift to $\GL_{2n}(\A)$ of the form
\begin{equation}
\label{E:weak-lift2}
 \BIGboxplus_i \phi_i.
\end{equation}
By the strong multiplicity one theorem \cite[Theorem 4.4]{js}, the two expressions \eqref{E:weak-lift1} and \eqref{E:weak-lift2} must agree.
However, $\tau_i$ in \eqref{E:weak-lift1} either is non-self-dual or is self-dual but occurs with multiplicity at least $2$, whereas $\phi_i$ in \eqref{E:weak-lift2} is self-dual and occurs with multiplicity $1$.
This is a contradiction.
Hence we have $\mathcal{V} \subset \AA_\cusp(\Mp_{2n})$ as required.
\end{proof}

\subsection{The multiplicity preservation}

As a consequence of Proposition \ref{p:tempered}, we deduce:

\begin{cor}
\label{c:mult-preserve}
Let $\phi$ be a generic elliptic $A$-parameter for $\Mp_{2n}$.
Suppose that
\[
 L^2_{\phi,\psi}(\Mp_{2n}) \cong \bigoplus_{\pi}  m_\pi \pi.
\]
If $r>2n+1$ and $\phi' = \phi \oplus S_{2r-2n}$, then
\[
 L^2_{\phi'}(\SO_{2r+1}) \cong \bigoplus_{\pi} m_\pi \theta_\psi^{\abs}(\pi).
\]
\end{cor}

\begin{proof}
For any irreducible genuine unitary representation $\pi$ of $\Mp_{2n}(\A)$ such that the $L$-parameter of $\pi_v$ (relative to $\psi_v$) is $\phi_v$ for almost all $v$, we have
\begin{equation}
\label{eq:mult-pre}
 m_\disc(\pi) = m_\disc(\theta_\psi^{\abs}(\pi))
\end{equation}
by Theorem \ref{T:inequalities} and Proposition \ref{p:tempered}.
In view of Remark \ref{rem:embed}, it remains to show that for any irreducible summand $\sigma$ of $L^2_{\phi'}(\SO_{2r+1})$, there exists an irreducible summand $\pi$ of $L^2_{\phi,\psi}(\Mp_{2n})$ such that $\sigma \cong \theta_\psi^{\abs}(\pi)$.
Since the $L$-parameter of $\sigma_v$ is
\[
 \varphi_{\phi'_v} = \phi_v \oplus 
 \Big( \bigoplus_{j=1}^{2r-2n} |\cdot|^{r-n+\frac{1}{2}-j} \Big)
\]
for almost all $v$, it follows from the local theta correspondence for unramified representations that
\begin{equation}
\label{eq:sigma_v}
 \sigma_v \cong \theta_{\psi_v}(\pi_{\phi_v}) 
\end{equation}
for almost all $v$, where $\pi_{\phi_v}$ is the irreducible genuine unramified representation of $\Mp_{2n}(F_v)$ whose $L$-parameter (relative to $\psi_v$) is $\phi_v$.
Note that such $\pi_{\phi_v}$ is unitary (see \cite{tadic} and Remark \ref{r:unitarity} below).
As explained in \S \ref{ss:low-rank}, such $\sigma_v$ is of rank $2n$, and hence so is $\sigma$.
Hence there exist a unique irreducible genuine unitary representation $\pi$ of $\Mp_{2n}(\A)$ and a unique quadratic automorphic character $\chi$ of $\A^\times$ such that
\[
 \sigma \cong \theta_\psi^\abs(\pi) \otimes (\chi \circ \nu).
\]
Then we have
\[
 \sigma_v \cong \theta_{\psi_v}(\pi_v) \otimes (\chi_v \circ \nu)
\]
for all $v$.
Recalling the bijection \eqref{eq:li-bij}, we deduce from this and \eqref{eq:sigma_v} that $\pi_v \cong \pi_{\phi_v}$ and $\chi_v$ is trivial for almost all $v$.
Since $\chi$ is automorphic, it must be trivial, so that $\sigma \cong \theta_\psi^{\abs}(\pi)$.
Hence, by \eqref{eq:mult-pre}, we have $m_\disc(\pi) = m_\disc(\sigma) > 0$.
\end{proof}

\subsection{Multiplicity formula for $\SO_{2r+1}$}

Corollary \ref{c:mult-preserve} allows us to infer a local-global structure of $L^2_{\phi,\psi}(\Mp_{2n})$ with a multiplicity formula from Arthur's result \cite{a}, which we now recall.
Let $\phi'$ be an elliptic $A$-parameter for $\SO_{2r+1}$ and write
\[
 \phi' = \bigoplus_i \phi_i \boxtimes S_{d_i} 
\]
as in \eqref{eq:A-param}.
Let
\[
 S_{\phi'} = \bigoplus_i (\Z/2 \Z) a'_i
\]
be a free $\Z/2\Z$-module with a basis $\{ a_i' \}$, where $a_i'$ corresponds to $\phi_i \boxtimes S_{d_i}$, and put 
\[
 \bar{S}_{\phi'} = S_{\phi'} / \Delta(\Z/2 \Z).
\]
For each place $v$ of $F$, we regard the localization $\phi'_v$ of $\phi'$ at $v$ as a local $A$-parameter $\phi'_v: L_{F_v} \times \SL_2(\C) \rightarrow \Sp_{2r}(\C)$.
Let $S_{\phi'_v}$ be the component group of the centralizer of the image of $\phi'_v$ in $\Sp_{2r}(\C)$ and put $\bar{S}_{\phi'_v} = S_{\phi'_v} / \langle z_{\phi'_v} \rangle$, where $z_{\phi'_v}$ is the image of $-1 \in \Sp_{2r}(\C)$ in $S_{\phi'_v}$.
Then we have a canonical map $\bar{S}_{\phi'} \rightarrow \bar{S}_{\phi'_v}$.
Thus, we obtain a compact group $\bar{S}_{\phi', \A} = \prod_v \bar{S}_{\phi'_v}$ equipped with the diagonal map $\Delta: \bar{S}_{\phi'} \rightarrow \bar{S}_{\phi', \A}$.

To each local $A$-parameter $\phi'_v$, Arthur \cite{a} assigned a finite set of (possibly zero, possibly reducible) semisimple representations of $\SO_{2r+1}(F_v)$ of finite length:
\[
 \Pi_{\phi'_v}(\SO_{2r+1}) = \{ \sigma_{\eta'_v} \, | \, \eta'_v \in \hat{\bar{S}}_{\phi'_v} \}.
\]
If $\phi'$ is generic, then $\Pi_{\phi'_v}(\SO_{2r+1})$ is simply the local $L$-packet associated to $\phi'_v$ (regarded as a local $L$-parameter) by the local Langlands correspondence.
In particular, if further $v$ is nonarchimedean or complex, then $\sigma_{\eta'_v}$ is nonzero and irreducible for any $\eta'_v$.
However, if $\phi'$ is not generic, then Arthur's result \cite{a} does not provide explicit knowledge of the representation $\sigma_{\eta'_v}$.
Fortunately, in the nonarchimedean case, M{\oe}glin's results \cite{m1,m2,m3,m4} provide an alternative explicit construction of $\Pi_{\phi'_v}(\SO_{2r+1})$ and rather precise knowledge of the properties of $\sigma_{\eta'_v}$ such as nonvanishing, multiplicity-freeness, and irreducibility.
These results will be reviewed in \S \ref{s:local1} below.

For any $\eta' = \bigotimes_v \eta'_v \in \hat{\bar{S}}_{\phi', \A}$, we may form a semisimple representation $\sigma_{\eta'} = \bigotimes_v \sigma_{\eta'_v}$ of $\SO_{2r+1}(\A)$.
Let $\epsilon_{\phi'}$ be the quadratic character of $\bar{S}_{\phi'}$ defined by \cite[(1.5.6)]{a}.
Then Arthur's multiplicity formula \cite[Theorem 1.5.2]{a} asserts that
\begin{equation}
\label{eq:mult-so}
 L^2_{\phi'}(\SO_{2r+1}) \cong \bigoplus_{\eta' \in \hat{\bar{S}}_{\phi', \A}} m_{\eta'} \sigma_{\eta'}, 
\end{equation}
where
\[
 m_{\eta'} =
 \begin{cases}
  1 & \text{if $\Delta^* \eta' = \epsilon_{\phi'}$;} \\
  0 & \text{otherwise.}
 \end{cases}
\]

\subsection{Structure of $L^2_{\phi,\psi}(\Mp_{2n})$}
\label{ss:near-eq-mp}

Finally, with the help of Corollary \ref{c:mult-preserve}, we can transfer the structure of $L^2_{\phi'}(\SO_{2r+1})$ to $L^2_{\phi,\psi}(\Mp_{2n})$ for any generic elliptic $A$-parameter $\phi$ for $\Mp_{2n}$, where $\phi' = \phi \oplus S_{2r-2n}$ with $r > 2n+1$.
If we write $\phi = \bigoplus_i \phi_i$ for some pairwise distinct irreducible symplectic cuspidal automorphic representations $\phi_i$ of $\GL_{n_i}(\A)$, then $S_\phi$ and $S_{\phi'}$ are of the form
\[
 S_\phi = \bigoplus_i (\Z/2 \Z) a_i, \qquad 
 S_{\phi'} = S_\phi  \oplus (\Z/2\Z) a'_0,
\]
where $a_i$ and $a'_0$ correspond to $\phi_i$ and $S_{2r-2n}$, respectively.
In particular, the natural map
\[
 \iota:S_\phi \longhookrightarrow S_{\phi'} \longtwoheadrightarrow \bar{S}_{\phi'}
\]
is an isomorphism.
Put $\epsilon_\phi = \iota^* \epsilon_{\phi'}$.

\begin{lem}
\label{l:epsilon}
We have
\[
 \epsilon_{\phi}(a_i) = \epsilon(\tfrac{1}{2}, \phi_i).
\]
Moreover, if the central $L$-value $L(\frac{1}{2}, \phi)$ does not vanish, then $\epsilon_{\phi}$ is trivial.
\end{lem}

\begin{proof}
Recall that $\phi_i$ is symplectic and hence $n_i$ is even (see \cite[\S 9]{js-aa}).
Let $\mathcal{L}_\phi = \prod_i \Sp_{n_i}(\C)$ be a substitute for the hypothetical Langlands group of $F$ as in \cite[(1.4.4)]{a}, so that for any place $v$ of $F$, we may regard $\phi_v$ as a local $L$-parameter $\phi_v : L_{F_v} \rightarrow \mathcal{L}_\phi$ via the local Langlands correspondence.
Let $\tilde{\phi} : \mathcal{L}_{\phi} \rightarrow \Sp_{2n}(\C)$ be a natural embedding.
We define a homomorphism $\tilde{\phi}' : \mathcal{L}_\phi \times \SL_2(\C) \rightarrow \Sp_{2r}(\C)$ by 
\[
 \tilde{\phi}' = \tilde{\phi} \oplus S_{2r-2n}, 
\]
where we regard $\tilde{\phi}$ and $\tilde{\phi}'$ as representations of $\mathcal{L}_\phi$ and $\mathcal{L}_\phi \times \SL_2(\C)$, respectively.
Then we have a natural embedding $\kappa:S_{\phi'} \hookrightarrow \Sp_{2r}(\C)$ whose image is equal to the centralizer of the image of $\tilde{\phi}'$ in $\Sp_{2r}(\C)$.
More explicitly, we have
\[
 \Std \circ (\kappa \times \tilde{\phi}') = \Big( \bigoplus_i \kappa_i \boxtimes \Std_i \Big) \oplus (\kappa_0' \boxtimes S_{2r-2n})
\]
as representations of $S_{\phi'} \times \mathcal{L}_\phi \times \SL_2(\C)$, where $\Std$ is the standard representation of $\Sp_{2r}(\C)$, $\Std_i$ is the composition of the projection $\mathcal{L}_\phi \rightarrow \Sp_{n_i}(\C)$ with the standard representation of $\Sp_{n_i}(\C)$, and $\kappa_i$ and $\kappa_0'$ are the characters of $S_{\phi'}$ given by
\begin{align*}
 \kappa_i(a_j) & =
 \begin{cases}
  -1 & \text{if $i = j$;} \\
  1 & \text{if $i \ne j$,}
 \end{cases} & 
 \kappa_i(a_0') & = 1, \\
 \kappa'_0(a_j) & = 1, &
 \kappa'_0(a_0') & = -1.
\end{align*}
We consider the orthogonal representation $\Ad \circ (\kappa \times \tilde{\phi}')$ of $S_{\phi'} \times \mathcal{L}_\phi \times \SL_2(\C)$, where $\Ad$ is the adjoint representation of $\Sp_{2r}(\C)$, and decompose it as
\[
 \Ad \circ (\kappa \times \tilde{\phi}') = \bigoplus_j \epsilon_j \boxtimes \rho_j \boxtimes S_{d_j}
\]
for some characters $\epsilon_j$ of $S_{\phi'}$, some irreducible representations $\rho_j$ of $\mathcal{L}_\phi$, and some positive integers $d_j$.
In fact, $\epsilon_j$ descends to a character of $\bar{S}_{\phi'}$ for all $j$.
If $\rho_j$ is symplectic, we denote by
\[
 \epsilon(s, \phi, \rho_j) = \prod_v \epsilon(s, \rho_j \circ \phi_v, \psi_v)
\]
the $\epsilon$-function of $\phi$ relative to $\rho_j$, which is independent of the choice of $\psi$.
Note that $\epsilon(\frac{1}{2}, \phi, \rho_j) \in \{ \pm 1 \}$.
Then $\epsilon_{\phi'}$ is given by
\[
 \epsilon_{\phi'} = \prod_j \epsilon_j,
\]
where the product runs over indices $j$ such that $\rho_j$ is symplectic and $\epsilon(\frac{1}{2}, \phi, \rho_j) = -1$.
On the other hand, we have
\[
 \Ad \circ \tilde{\phi}' = (\Ad \circ \tilde{\phi}) \oplus (\tilde{\phi} \boxtimes S_{2r-2n})
 \oplus S_3 \oplus S_7 \oplus \cdots \oplus S_{4r-4n-1}
\]
as representations of $\mathcal{L}_\phi \times \SL_2(\C)$.
Since $\Ad \circ \tilde{\phi}$ is orthogonal, it does not contribute to $\epsilon_{\phi'}$.
Moreover, we have
\[
 \tilde{\phi} \boxtimes S_{2r-2n} = \bigoplus_i \kappa_i \kappa_0' \boxtimes \Std_i \boxtimes S_{2r-2n}
\]
as representations of $S_{\phi'} \times \mathcal{L}_\phi \times \SL_2(\C)$.
Since $\Std_i$ is symplectic and $\epsilon(s, \phi, \Std_i) = \epsilon(s, \phi_i)$, this implies the first assertion.

If $L(\frac{1}{2}, \phi) \ne 0$, then since $L(\frac{1}{2}, \phi) = \prod_i L(\frac{1}{2}, \phi_i)$, we have $L(\frac{1}{2}, \phi_i) \ne 0$ and hence $\epsilon(\tfrac{1}{2}, \phi_i) = 1$ for all $i$.
Thus, the second assertion follows from the first one.
\end{proof}

Similarly, for any place $v$ of $F$, the natural map
\[
 \iota_v:S_{\phi_v} \longhookrightarrow S_{\phi'_v} \longtwoheadrightarrow \bar{S}_{\phi'_v}
\]
is an isomorphism.
For any $\eta_v' \in \hat{\bar{S}}_{\phi'_v}$, we write 
\[
 \sigma_{\eta'_v} = \bigoplus_i m_{\eta'_v, i} \sigma_{\eta'_v, i}
\]
for some positive integers $m_{\eta'_v, i}$ and some pairwise distinct irreducible representations $\sigma_{\eta'_v, i}$ of $\SO_{2r+1}(F_v)$.
Put $\eta_v = \iota^*_v \eta_v'$ and 
\[
 \tilde\pi_{\eta_v} = \bigoplus_i m_{\eta'_v, i} \theta_{\psi_v}(\sigma_{\eta'_v, i}),
\]
where $\theta_{\psi_v}(\sigma_{\eta'_v, i})$ is the theta lift to $\Mp_{2n}(F_v)$.
Thus, for any $\eta = \bigotimes_v \eta_v \in \hat{S}_{\phi, \A}$, we may form a semisimple genuine representation $\tilde\pi_{\eta} = \bigotimes_v \tilde\pi_{\eta_v}$ of $\Mp_{2n}(\A)$.

\begin{prop}
\label{p:mult-mp}
Let $\phi$ be a generic elliptic $A$-parameter for $\Mp_{2n}$.
Then we have
\[
 L^2_{\phi,\psi}(\Mp_{2n}) \cong \bigoplus_{\eta \in \hat{S}_{\phi,\A}} m_\eta \tilde\pi_{\eta},
\]
where
\[
 m_\eta = 
 \begin{cases}
  1 & \text{if $\Delta^* \eta = \epsilon_{\phi}$;} \\
  0 & \text{otherwise.}
 \end{cases}
\]
\end{prop}

\begin{proof}
Suppose that
\[
 L^2_{\phi'}(\SO_{2r+1}) \cong \bigoplus_\sigma m_\sigma \sigma.
\]
Then, by Corollary \ref{c:mult-preserve} and the Howe duality, we have
\[
 L^2_{\phi,\psi}(\Mp_{2n}) \cong \bigoplus_\sigma m_\sigma \theta_\psi^\abs(\sigma).
\]
This and the multiplicity formula \eqref{eq:mult-so} imply the assertion.
\end{proof}

Hence, to complete the proof of Theorem \ref{t:main2},
it remains to describe $\tilde\pi_\eta$ in terms of the local Shimura correspondence.
This will be established in Proposition \ref{p:key} below.

\begin{rem}
In the above argument, we have fixed an integer $r > 2n+1$ and do not know a priori that $\tilde\pi_\eta$ is independent of the choice of $r$.
This seems not immediate and will be deduced from the description of $\tilde\pi_\eta$ in terms of the local Shimura correspondence.
While this is a purely local problem, we will address it by a global argument in \S \ref{s:local2} below.
\end{rem}

\section{\textbf{Local $L$- and $A$-packets}}
\label{s:local1}

In this section, we review the representation theory of metaplectic and orthogonal groups over local fields.
In particular, we state irreducibility of some induced representations which will play an important role in an inductive argument in \S \ref{s:local2} below.

\subsection{$L$-parameters}

Let $F$ be a local field of characteristic zero and put 
\[
 L_F =
 \begin{cases}
  \text{the Weil group of $F$} & \text{if $F$ is archimedean;} \\
  \text{the Weil--Deligne group of $F$} & \text{if $F$ is nonarchimedean.}
 \end{cases}
\]
Then the local Langlands correspondence \cite{langlands2,ht,henniart,scholze} provides a bijection
\[
 \Irr \GL_n \longleftrightarrow 
 \{ \text{$n$-dimensional representations of $L_F$} \}.
\]
Let $\phi$ be an $n$-dimensional representation of $L_F$.
We may regard $\phi$ as an $L$-parameter $\phi:L_F \rightarrow \GL_n(\C)$.
We say that:
\begin{itemize}
 \item $\phi$ is symplectic if there exists a nondegenerate antisymmetric bilinear form $b:\C^n \times \C^n \rightarrow \C$ such that $b(\phi(w)x, \phi(w)y) = b(x,y)$ for all $w \in L_F$ and $x, y \in \C^n$;
 \item $\phi$ is tempered if the image of the Weil group of $F$ under $\phi$ is relatively compact in $\GL_n(\C)$.
\end{itemize}
If $\phi$ is irreducible and symplectic, then it is tempered.
Let $\tau$ be the irreducible representation of $\GL_n$ associated to $\phi$.
Then we have:
\begin{itemize}
 \item $\tau$ is essentially square-integrable if and only if $\phi$ is irreducible;
 \item $\tau$ is tempered if and only if $\phi$ is tempered.
\end{itemize}

Let $\phi$ be a $2n$-dimensional symplectic representation of $L_F$.
We may regard $\phi$ as an $L$-parameter $\phi:L_F \rightarrow \Sp_{2n}(\C)$.
We decompose $\phi$ as
\begin{equation}
\label{eq:L-param}
 \phi = \bigoplus_i m_i \phi_i 
\end{equation}
for some positive integers $m_i$ and some pairwise distinct irreducible representations $\phi_i$ of $L_F$.
We say that:
\begin{itemize}
 \item $\phi$ is good if $\phi_i$ is symplectic for all $i$;
 \item $\phi$ is tempered if $\phi_i$ is tempered for all $i$;
 \item $\phi$ is almost tempered if $\phi_i|\cdot|^{-s_i}$ is tempered for some $s_i \in \R$ with $|s_i| <\frac{1}{2}$ for all $i$.
\end{itemize}
If $\phi$ is good, then it is tempered.
Also, any localization of a global generic elliptic $A$-parameter for $\Mp_{2n}$ is almost tempered.

For any $2n$-dimensional symplectic representation $\phi$ of $L_F$, we may write
\begin{equation}
\label{eq:psi-decomp}
 \phi = \varphi \oplus \varphi^\vee \oplus \phi_0,
\end{equation}
where
\begin{itemize}
 \item $\varphi$ is a $k$-dimensional representation of $L_F$ whose irreducible summands are all non-symplectic;
 \item $\phi_0$ is a $2 n_0$-dimensional good symplectic representation of $L_F$;
 \item $k+n_0=n$.
\end{itemize}
More explicitly, if $\phi$ is of the form \eqref{eq:L-param}, then $\phi_0$ is given by
\[
 \phi_0 = \bigoplus_{i \in I_0} m_i \phi_i,
\]
where $I_0 = \{ i \, | \, \text{$\phi_i$ is symplectic} \}$.
Let $S_\phi$ be the component group of the centralizer of the image of $\phi$ in $\Sp_{2n}(\C)$ and $z_\phi$ the image of $-1 \in \Sp_{2n}(\C)$ in $S_\phi$.
Then $S_\phi$ is a free $\Z/2\Z$-module of the form
\[
 S_\phi = \bigoplus_{i \in I_0} (\Z/2\Z) a_i,
\]
where $a_i$ corresponds to $\phi_i$, and $z_\phi = (m_i a_i)$.
In particular, we have a natural identification $S_{\phi_0} = S_\phi$.

\subsection{Representation theory of $\SO(V)$}
\label{ss:rep-so}

The local Langlands correspondence \cite{langlands2,shelstad1,shelstad2,a,mr4} provides a partition
\[
 \Irr \SO(V) = \bigsqcup_\phi \Pi_\phi(\SO(V))
\]
and bijections
\begin{equation}
\label{eq:bij-so}
 \bigsqcup_V \Pi_\phi(\SO(V)) \longleftrightarrow \hat{S}_\phi, 
\end{equation}
where the first disjoint union runs over equivalence classes of $2n$-dimensional symplectic representations $\phi$ of $L_F$ and the second disjoint union runs over isometry classes of $(2n+1)$-dimensional quadratic spaces $V$ over $F$ with trivial discriminant.
For any $\sigma \in \Pi_\phi(\SO(V))$, we have:
\begin{itemize}
 \item $\sigma$ is square-integrable if and only if $\phi$ is good and multiplicity-free;
 \item $\sigma$ is tempered if and only if $\phi$ is tempered.
\end{itemize}
We write $\sigma = \sigma_\eta$ if $\sigma$ corresponds to $\eta \in \hat{S}_\phi$ under the bijection \eqref{eq:bij-so}.
Let $\hat{S}_{\phi,V} \subset \hat{S}_\phi$ be the subset of all $\eta$ such that $\sigma_\eta$ is a representation of $\SO(V)$.
Then we have
\[
 \hat{S}_{\phi,V} \subset \{ \eta \in \hat{S}_\phi \, | \, \eta(z_\phi) = \varepsilon(V) \},
\]
with equality when $F$ is nonarchimedean or $F=\C$.
Note that if $F = \C$, then $S_\phi = \{ 0 \}$, so that the condition $\eta(z_\phi) = \varepsilon(V)$ disappears.

Write $\phi = \varphi \oplus \varphi^\vee \oplus \phi_0$ as in \eqref{eq:psi-decomp}.
We have $\Pi_\phi(\SO(V)) = \varnothing$ unless the $F$-rank of $\SO(V)$ is greater than or equal to $k$, in which case there exists a $(2n_0+1)$-dimensional quadratic space $V_0$ over $F$ with trivial discriminant such that $V = \mathbb{H}^k \oplus V_0$.
Let $Q$ be a parabolic subgroup of $\SO(V)$ with Levi component $\GL_k \times \SO(V_0)$.
Let $\tau$ be the irreducible representation of $\GL_k$ associated to $\varphi$.
Then, by the inductive property of the local Langlands correspondence, for any $\eta \in \hat{S}_{\phi,V}$, $\sigma_\eta$ is equal to an irreducible subquotient of
\[
 \Ind^{\SO(V)}_Q(\tau \otimes \sigma_{\eta_0}),
\]
where we write $\eta$ as $\eta_0$ if we regard it as a character of $S_{\phi_0}$ via the identification $S_{\phi_0} = S_\phi$.

\begin{lem}
\label{l:irred-so}
If $\phi$ is almost tempered, then $\Ind^{\SO(V)}_Q(\tau \otimes \sigma_0)$ is irreducible for any $\sigma_0 \in \Pi_{\phi_0}(\SO(V_0))$.
\end{lem}

\begin{proof}
If $\phi$ is tempered, then by definition, $\Pi_\phi(\SO(V))$ consists of all irreducible summands of $\Ind^{\SO(V)}_Q(\tau \otimes \sigma_0)$ for all $\sigma_0 \in \Pi_{\phi_0}(\SO(V_0))$.
Hence the irreducibility of $\Ind^{\SO(V)}_Q(\tau \otimes \sigma_0)$ is a consequence of the bijectivity of \eqref{eq:bij-so} and the fact that $S_{\phi_0} = S_\phi$.

Thus, if $\phi$ is almost tempered, then by induction in stages, it remains to show that $\Ind^{\SO(V)}_Q(\tau \otimes \sigma_0)$ is irreducible if
\begin{itemize}
 \item $\tau$ is an irreducible representation of $\GL_k$ whose $L$-parameter is of the form $\bigoplus_i \phi_i |\cdot|^{s_i}$ for some tempered representations $\phi_i$ of $L_F$ and some $s_i \in \R$ with $0 < |s_i| < \frac{1}{2}$;
 \item $\sigma_0$ is an irreducible tempered representation of $\SO(V_0)$.
\end{itemize}
If $F$ is archimedean, then this follows from a result of Speh--Vogan \cite{spehv} (see also \cite[Chapter 8]{vogan81}).
If $F$ is nonarchimedean, then this follows from a result of M{\oe}glin--Waldspurger \cite[\S 2.14]{mw} and a conjecture of Gross--Prasad and Rallis \cite[Conjecture 2.6]{gp}, which is proved in \cite[Appendix B]{gi2}.
\end{proof}

\subsection{Representation theory of $\Mp_{2n}$}
\label{ss:rep-mp}

Fix a nontrivial additive character $\psi$ of $F$.
The local Shimura correspondence \cite{ab1,ab2,gs} asserts that the theta lift (relative to $\psi$) induces a bijection
\[
 \theta_\psi : \Irr \Mp_{2n} \longleftrightarrow \bigsqcup_V \Irr \SO(V)
\]
satisfying natural properties, where the disjoint union runs over isometry classes of $(2n+1)$-dimensional quadratic spaces $V$ over $F$ with trivial discriminant.
Composing this with the local Langlands correspondence for $\SO(V)$, we obtain a partition
\[
 \Irr \Mp_{2n} = \bigsqcup_\phi \Pi_{\phi,\psi}(\Mp_{2n})
\]
and bijections
\begin{equation}
\label{eq:bij-mp}
 \Pi_{\phi,\psi}(\Mp_{2n}) \longleftrightarrow \hat{S}_\phi,
\end{equation}
where the disjoint union runs over equivalence classes of $2n$-dimensional symplectic representations $\phi$ of $L_F$.
Since $\theta_\psi$ preserves the square-integrability and the temperedness of representations, for any $\pi \in \Pi_{\phi,\psi}(\Mp_{2n})$, we have:
\begin{itemize}
 \item $\pi$ is square-integrable if and only if $\phi$ is good and multiplicity-free;
 \item $\pi$ is tempered if and only if $\phi$ is tempered.
\end{itemize}
We write $\pi = \pi_\eta$ if $\pi$ corresponds to $\eta \in \hat{S}_\phi$ under the bijection \eqref{eq:bij-mp}.

Write $\phi = \varphi \oplus \varphi^\vee \oplus \phi_0$ as in \eqref{eq:psi-decomp}.
Let $P$ be a parabolic subgroup of $\Sp_{2n}$ with Levi component $\GL_k \times \Sp_{2n_0}$ and $\tilde{P}$ the preimage of $P$ in $\Mp_{2n}$.
Let $\tau$ be the irreducible representation of $\GL_k$ associated to $\varphi$ and $\tilde{\tau}_\psi = \tau \otimes \chi_\psi$ its twist by the genuine quartic character $\chi_\psi$ of the two-fold cover of $\GL_k$ given in \cite[\S 2.6]{gi1}.
Then, by the inductive property of the local Shimura correspondence, for any $\eta \in \hat{S}_\phi$, $\pi_\eta$ is equal to an irreducible subquotient of
\[
 \Ind^{\Mp_{2n}}_{\tilde{P}}(\tilde{\tau}_\psi \otimes \pi_{\eta_0}),
\]
where we write $\eta$ as $\eta_0$ if we regard it as a character of $S_{\phi_0}$ via the identification $S_{\phi_0} = S_\phi$.

\begin{lem}
\label{l:irred-mp}
If $\phi$ is almost tempered, then $\Ind^{\Mp_{2n}}_{\tilde{P}}(\tilde{\tau}_\psi \otimes \pi_0)$ is irreducible for any $\pi_0 \in \Pi_{\phi_0,\psi}(\Mp_{2n_0})$.
\end{lem}

\begin{proof}
If $\phi$ is tempered, then by the inductive property of the local Shimura correspondence, $\Pi_{\phi,\psi}(\Mp_{2n})$ consists of all irreducible summands of $\Ind^{\Mp_{2n}}_{\tilde{P}}(\tilde{\tau}_\psi \otimes \pi_0)$ for all $\pi_0 \in \Pi_{\phi_0,\psi}(\Mp_{2n_0})$.
Hence the irreducibility of $\Ind^{\Mp_{2n}}_{\tilde{P}}(\tilde{\tau}_\psi \otimes \pi_0)$ is a consequence of the bijectivity of \eqref{eq:bij-mp} and the fact that $S_{\phi_0} = S_\phi$.

Thus, if $\phi$ is almost tempered, then by induction in stages, it remains to show that $\Ind^{\Mp_{2n}}_{\tilde{P}}(\tilde{\tau}_\psi \otimes \pi_0)$ is irreducible if
\begin{itemize}
 \item $\tau$ is an irreducible representation of $\GL_k$ whose $L$-parameter is of the form $\bigoplus_i \phi_i |\cdot|^{s_i}$ for some tempered representations $\phi_i$ of $L_F$ and some $s_i \in \R$ with $0 < |s_i| < \frac{1}{2}$;
 \item $\pi_0$ is an irreducible tempered representation of $\Mp_{2n_0}$.
\end{itemize}
If $F = \C$, then this follows from a result of Speh--Vogan \cite{spehv} (see also \cite[Chapter 8]{vogan81}).
If $F = \R$, then the argument in \cite{spehv} should also work, but they only consider real reductive \emph{linear} Lie groups.
For the sake of completeness, we give a proof in \cite{gi-mp-real}.
Suppose that $F$ is nonarchimedean.
Then this follows from a result of Atobe \cite[Theorem 3.13]{atobe}, but we include the proof for the convenience of the reader.
We may assume that $\pi = \Ind^{\Mp_{2n}}_{\tilde{P}}(\tilde{\tau}_\psi \otimes \pi_0)$ is a standard module.
Let $\pi'$ be the unique irreducible quotient of $\pi$.
By the local Shimura correspondence, the theta lift $\sigma_0 = \theta_\psi(\pi_0)$ to $\SO(V_0)$ is nonzero and tempered for a unique $(2n_0+1)$-dimensional quadratic space $V_0$ over $F$ with trivial discriminant.
Put $\sigma = \Ind^{\SO(V)}_Q(\tau \otimes \sigma_0)$, where $V = \mathbb{H}^k \oplus V_0$ and $Q$ is the standard parabolic subgroup of $\SO(V)$ with Levi component $\GL_k \times \SO(V_0)$.
As shown in the proof of Lemma \ref{l:irred-so}, $\sigma$ is irreducible.
Then, by \cite{gs}, we have
\[
 \theta_\psi(\pi') = \sigma.  
\]
On the other hand, as in the proof of \cite[Theorem 8.1]{gs}, it follows from Kudla's filtration (see e.g.~\cite[Proposition 7.3]{gs}) and \cite[Lemma 7.4]{gs}, which continues to hold for $s>-\frac{1}{2}$, that there exists a surjection
\[
 \Ind^{\Mp_{2n}}_{\tilde{P}}((\tau^\vee \otimes \chi_\psi) \otimes \Theta_\psi(\sigma_0)) \longtwoheadrightarrow \Theta_\psi(\sigma).
\]
Since $\Theta_\psi(\sigma_0) = \theta_\psi(\sigma_0) = \pi_0$ by \cite{gs} and $\theta_\psi(\sigma) = \pi'$, this induces a surjection
\[
 \Ind^{\Mp_{2n}}_{\tilde{P}}((\tau^\vee \otimes \chi_\psi) \otimes \pi_0) \longtwoheadrightarrow \pi'.
\]
Taking the contragredient and applying the MVW involution \cite{mvw,sun}, we obtain an injection 
\[
 \pi' \longhookrightarrow \Ind^{\Mp_{2n}}_{\tilde{P}}(\tilde{\tau}_\psi \otimes \pi_0).
\]
Hence $\pi'$ occurs in $\pi$ as a quotient and as a subrepresentation.
Since $\pi'$ occurs in $\pi$ with multiplicity one (see e.g.~\cite[Remark 4.5]{bj}), $\pi$ must be irreducible.
\end{proof}

\begin{rem}
\label{r:unitarity}
Suppose that $F$ is nonarchimedean of odd residual characteristic and $\psi$ is of order zero.
Then we may regard $K = \Sp_{2n}(\mathcal{O})$ as a subgroup of $\Mp_{2n}$ via the standard splitting, where $\mathcal{O}$ is the integer ring of $F$.
If an irreducible genuine representation $\pi$ of $\Mp_{2n}$ is unramified, i.e.~has a nonzero $K$-fixed vector, then it must be the unique unramified subquotient of a principal series representation
\begin{equation}
\label{eq:unram-ps}
 \Ind^{\Mp_{2n}}_{\tilde{B}}(\chi_{\psi} |\cdot|^{s_1} \otimes \cdots \otimes \chi_{\psi} |\cdot|^{s_n}), 
\end{equation}
where $B$ is a Borel subgroup of $\Sp_{2n}$, $\tilde{B}$ is the preimage of $B$ in $\Mp_{2n}$, and $s_1, \ldots, s_n \in \C$.
In this case, by the local theta correspondence for unramified representations, the $L$-parameter of $\pi$ (relative to $\psi$) is
\[
 \bigoplus_{i=1}^n (|\cdot|^{s_i} \oplus |\cdot|^{-s_i}).
\]
Also, as shown in the proof of Lemma \ref{l:irred-mp}, the representation \eqref{eq:unram-ps} is irreducible if $|\Re s_i|<\frac{1}{2}$ for all $i$.
Hence, by \cite[Lemma 3.3]{lmt}, the representation \eqref{eq:unram-ps} is unitary if either $s_i \in \sqrt{-1} \R$ or $s_i \in \R$ with $|s_i|<\frac{1}{2}$ for all $i$.
\end{rem}

\subsection{Some $A$-packets for $\SO_{2r+1}$}

Let $\phi'$ be a $2r$-dimensional symplectic representation of $L_F \times \SL_2(\C)$.
We may regard $\phi'$ as an $A$-parameter $\phi':L_F \times \SL_2(\C) \rightarrow \Sp_{2r}(\C)$ and associate to it an $L$-parameter $\varphi_{\phi'} : L_F \rightarrow \Sp_{2r}(\C)$ by
\begin{equation}
\label{eq:assoc-Lparam}
 \varphi_{\phi'}(w) = \phi' \!
 \left( w,
 \begin{pmatrix}
  |w|^{\frac{1}{2}} & \\
  & |w|^{-\frac{1}{2}}
 \end{pmatrix}
 \right).
\end{equation}
Let $S_{\phi'}$ be the component group of the centralizer of the image of $\phi'$ in $\Sp_{2r}(\C)$ and put $\bar{S}_{\phi'} = S_{\phi'} / \langle z_{\phi'} \rangle$, where $z_{\phi'}$ is the image of $-1 \in \Sp_{2r}(\C)$ in $S_{\phi'}$.
Then Arthur \cite{a} assigned to $\phi'$ an $A$-packet:
\[
 \Pi_{\phi'}(\SO_{2r+1}) = \{ \sigma_{\eta'} \, | \, \eta' \in \hat{\bar{S}}_{\phi'} \},
\]
where $\sigma_{\eta'}$ is a (possibly zero, possibly reducible) semisimple representation of $\SO_{2r+1}$ of finite length.

From now on, we only consider an $A$-parameter $\phi'$ of the form
\[
 \phi' = \phi \oplus S_{2r-2n}
\]
for some $2n$-dimensional symplectic representation $\phi$ of $L_F$ with $2n < r-1$.
Then we have:

\begin{prop}
\label{p:moeglin}
Assume that $\phi$ is almost tempered.
Then $\sigma_{\eta'}$ is nonzero and irreducible for any $\eta' \in \hat{\bar{S}}_{\phi'}$.
Moreover, $\Pi_{\phi'}(\SO_{2r+1})$ is multiplicity-free, i.e.~the $\sigma_{\eta'}$'s are pairwise distinct.
\end{prop}

This proposition is largely due to M{\oe}glin \cite{m1,m2,m3,m4} when $F$ is nonarchimedean, to M{\oe}glin \cite{m5} and M{\oe}glin--Renard \cite{mr} when $F=\C$, and to Arancibia--M{\oe}glin--Renard \cite{amr} and M{\oe}glin--Renard \cite{mr2,mr3} when $F = \R$ and $\phi$ is good.
The reader can also consult an efficient and concise exposition of M{\oe}glin's results by B.~Xu \cite{xu}.
In the next section, we will give a proof of Proposition \ref{p:moeglin} based on theta lifts.
For this, we need the following irreducibility of some induced representations.

Write $\phi = \varphi \oplus \varphi^\vee \oplus \phi_0$ as in \eqref{eq:psi-decomp}.
Put $\phi'_0 = \phi_0 \oplus S_{2r-2n}$, so that $\phi' = \varphi \oplus \varphi^\vee \oplus \phi_0'$.
Let $Q'$ be a parabolic subgroup of $\SO_{2r+1}$ with Levi component $\GL_k \times \SO_{2r-2k+1}$.
Let $\tau$ be the irreducible representation of $\GL_k$ associated to $\varphi$.
Then, by the definition and the inductive property of $A$-packets (see \cite{a}), for any $\eta' \in \hat{\bar{S}}_{\phi'}$, $\sigma_{\eta'}$ is equal to the semisimplification of
\[
 \Ind^{\SO_{2r+1}}_{Q'}(\tau \otimes \sigma_{\eta'_0}),
\]
where we write $\eta'$ as $\eta'_0$ if we regard it as a character of $\bar{S}_{\phi'_0}$ via the identification $\bar{S}_{\phi'_0} = \bar{S}_{\phi'}$.

\begin{lem}
\label{l:irred}
Assume that $\phi$ is almost tempered.
Then $\Ind^{\SO_{2r+1}}_{Q'}(\tau \otimes \sigma'_0)$ is irreducible for any irreducible subrepresentation $\sigma_0'$ of any representation in the $A$-packet $\Pi_{\phi_0'}(\SO_{2r-2k+1})$.
\end{lem}

\begin{proof}
The assertion was proved in a more general context by M{\oe}glin \cite[\S 3.2]{m3}, \cite[Proposition 5.1]{m4} when $F$ is nonarchimedean, and by M{\oe}glin--Renard \cite[\S 6]{mr} when $F = \C$.
In \cite{gi-mp-real}, we give a proof based on the Kazhdan--Lusztig algorithm when $F = \R$, noting that $\sigma_0'$ has half integral infinitesimal character by \cite[Lemme 3.4]{bergeron-clozel}.
\end{proof}

\section{\textbf{Comparison of local theta lifts}}
\label{s:local2}

As explained in \S \ref{ss:near-eq-mp}, to finish the proof of Theorem \ref{t:main2}, it remains to describe the local theta lift from $\SO_{2r+1}$ to $\Mp_{2n}$ with $r > 2n+1$ in terms of the local Shimura correspondence.
Namely, we need to compare the local theta correspondences for the following reductive dual pairs:
\begin{itemize}
 \item $(\Mp_{2n}, \SO_{2n+1})$ in the equal rank case (and its inner forms);
 \item $(\Mp_{2n}, \SO_{2r+1})$ in the stable range.
\end{itemize}
To distinguish them, we keep using $\theta_\psi$ to denote the theta correspondence for the former but change it to $\vartheta_\psi$ for the latter.
We emphasize that we need all inner forms of $\SO_{2n+1}$ in the former but only need the split form $\SO_{2r+1}$ in the latter.

Let $F$ be a local field of characteristic zero and fix a nontrivial additive character $\psi$ of $F$.
Let $\phi$ be a $2n$-dimensional symplectic representation of $L_F$ and put $\phi' = \phi \oplus S_{2r-2n}$ with $r > 2n+1$.
Then we have a natural isomorphism
\[
 \iota : S_\phi \longhookrightarrow S_{\phi'} \longtwoheadrightarrow \bar{S}_{\phi'}.
\]

We now state the main result of this section.

\begin{prop}
\label{p:key}
Assume that $\phi$ is almost tempered.
Then, for any $\eta' \in \hat{\bar{S}}_{\phi'}$, we have
\[
 \vartheta_\psi(\sigma_{\eta'}) = \pi_\eta,
\]
where $\eta = \iota^* \eta'$.
\end{prop}

The rest of this section is devoted to the proof of Propositions \ref{p:moeglin} and \ref{p:key}.

\subsection{Reduction to the case of good $L$-parameters}

We proceed by induction on $n$.

\begin{lem}
Propositions \ref{p:moeglin} and \ref{p:key} hold for $n=0$.
\end{lem}

\begin{proof}
If $n=0$, then $\phi' = S_{2r}$ and $\Pi_{\phi'}(\SO_{2r+1}) = \{ \sigma_\1 \}$, where $\sigma_\1$ is the trivial representation of $\SO_{2r+1}$.
Moreover, the theta lift $\vartheta_\psi(\sigma_\1)$ to $\Mp_0$ is the genuine $1$-dimensional representation of $\Mp_0 = \{ \pm 1 \}$, which is also the theta lift of the trivial representation of $\SO_1 = \{ 1 \}$.
This completes the proof.
\end{proof}

From now on, we assume that $n>0$.

\begin{lem}
\label{l:red-good}
Assume that Propositions \ref{p:moeglin} and \ref{p:key} hold for all $2n_0$-dimensional good symplectic representations of $L_F$ for all $n_0 < n$.
Then they also hold for all $2n$-dimensional almost tempered non-good symplectic representations of $L_F$.
\end{lem}

\begin{proof}
Let $\phi$ be a $2n$-dimensional almost tempered non-good symplectic representation of $L_F$.
Write $\phi = \varphi \oplus \varphi^\vee \oplus \phi_0$ as in \eqref{eq:psi-decomp}.
In particular, $\phi_0$ is a $2 n_0$-dimensional good symplectic representation of $L_F$ with $n_0 < n$.
Put $\phi'_0 = \phi_0 \oplus S_{2r-2n}$, so that $\phi' = \varphi \oplus \varphi^\vee \oplus \phi_0'$.
Let $\eta' \in \hat{\bar{S}}_{\phi'}$ and put $\eta = \iota^* \eta'$.
We write $\eta'$ (resp.~$\eta$) as $\eta'_0$ (resp.~$\eta_0$) if we regard it as a character of $\bar{S}_{\phi_0'}$ (resp.~$S_{\phi_0}$) via the natural identification.
Then, by definition, $\sigma_{\eta'}$ is the semisimplification of $\Ind^{\SO_{2r+1}}_{Q'}(\tau \otimes \sigma_{\eta'_0})$, where $Q'$ is the standard parabolic subgroup of $\SO_{2r+1}$ with Levi component $\GL_k \times \SO_{2r-2k+1}$ and $\tau$ is the irreducible representation of $\GL_k$ associated to $\varphi$.
Since $\sigma_{\eta'_0}$ is nonzero and irreducible by assumption, so is $\sigma_{\eta'}$ by Lemma \ref{l:irred}.
Moreover, the theta lift $\vartheta_\psi(\sigma_{\eta'_0})$ to $\Mp_{2n_0}$ is $\pi_{\eta_0}$ by assumption.
Hence it follows from the induction principle \cite{kudla}, \cite[Corollary 3.21]{ab1}, \cite[Theorem 8.4]{ab2}, which easily extends to the case at hand, that there exists a nonzero equivariant map
\[
 \omega_\psi \longrightarrow 
 \Ind^{\Mp_{2n}}_{\tilde{P}}(\tilde{\tau}_\psi \otimes \pi_{\eta_0})
 \otimes
 \Ind^{\SO_{2r+1}}_{Q'}(\tau \otimes \sigma_{\eta'_0}),
\]
where $P$ is the standard parabolic subgroup of $\Sp_{2n}$ with Levi component $\GL_k \times \Sp_{2n_0}$.
Since $\Ind^{\Mp_{2n}}_{\tilde{P}}(\tilde{\tau}_\psi \otimes \pi_{\eta_0})$ is irreducible by Lemma \ref{l:irred-mp}, this implies that the theta lift $\vartheta_\psi(\sigma_{\eta'})$ to $\Mp_{2n}$ is $\pi_{\eta}$.
Thus, since $\Pi_{\phi,\psi}(\Mp_{2n})$ is multiplicity-free, so is $\Pi_{\phi'}(\SO_{2r+1})$.
This completes the proof.
\end{proof}

We may now assume that Propositions \ref{p:moeglin} and \ref{p:key} hold for all $2n$-dimensional almost tempered non-good symplectic representations of $L_F$.
It remains to show that they also hold for all $2n$-dimensional good symplectic representations of $L_F$.
In particular, we have finished the proof for $F=\C$ since any irreducible representation of $L_\C$ is $1$-dimensional and hence non-symplectic.

Later, we also need the following description of the theta lift from $\SO_{2r+1}$ to $\Mp_{2n}$.

\begin{lem}
\label{l:theta-triv-char}
Assume that $F$ is nonarchimedean or $F=\C$, and that $\phi$ is almost tempered.
Let $\varphi_{\phi'}$ be the $L$-parameter associated to $\phi'$ by \eqref{eq:assoc-Lparam} and $\sigma'$ the irreducible representation in the $L$-packet $\Pi_{\varphi_{\phi'}}(\SO_{2r+1})$ associated to the trivial character of $S_{\varphi_{\phi'}}$.
Then we have
\[
 \vartheta_\psi(\sigma') = \pi_\1,
\]
where $\pi_\1$ is the irreducible representation in the $L$-packet $\Pi_{\phi,\psi}(\Mp_{2n})$ associated to the trivial character of $S_\phi$.
\end{lem}

\begin{proof}
If $\phi$ is non-good, then it follows from \cite[Proposition 7.4.1]{a} and Proposition \ref{p:moeglin} that $\sigma'$ is the representation in the $A$-packet $\Pi_{\phi'}(\SO_{2r+1})$ associated to the trivial character of $\bar{S}_{\phi'}$.
Hence, by Proposition \ref{p:key}, we have $\vartheta_\psi(\sigma') = \pi_\1$.

We may now assume that $\phi$ is good (and hence tempered), so that $F$ is nonarchimedean.
By definition, $\sigma'$ is the unique irreducible quotient of the standard module
\[
 \Ind^{\SO_{2r+1}}_{Q'}(|\cdot|^{r-n-\frac{1}{2}} \otimes |\cdot|^{r-n-\frac{3}{2}} \otimes \dots \otimes |\cdot|^{\frac{1}{2}} \otimes \sigma),
\]
where $Q'$ is the standard parabolic subgroup of $\SO_{2r+1}$ with Levi component $(\GL_1)^{r-n} \times \SO_{2n+1}$ and $\sigma$ is the irreducible tempered representation in the $L$-packet $\Pi_{\phi}(\SO_{2n+1})$ associated to the trivial character of $S_\phi$.
Since the theta lift of $\sigma$ to $\Mp_{2n}$ is $\pi_\1$, we have $\vartheta_\psi(\sigma') = \pi_\1$ by \cite[Proposition 3.2]{gt1}.
\end{proof}

\subsection{Multiplicity formula and globalization}
\label{ss:mult-global}

To finish the proof of Propositions \ref{p:moeglin} and \ref{p:key}, we appeal to a global argument.
One of the global ingredients is Arthur's multiplicity formula \eqref{eq:mult-so} and the following variant.

Let $\F$ be a number field and $\A$ the ad\`ele ring of $\F$.
Let $\V$ be a $(2n+1)$-dimensional quadratic space over $\F$ with trivial discriminant.
Recall the (expected) decomposition of $L^2_{\disc}(\SO(\V))$ into near equivalence classes described in \S \ref{ss:near-eq-so}.
We only consider the near equivalence class $L^2_{\varPhi}(\SO(\V))$ associated to a generic elliptic $A$-parameter $\varPhi$ for $\SO(\V)$.
As in \S \ref{ss:mult-mp}, we formally define the global component group $S_\varPhi$ of $\varPhi$ equipped with a canonical map $S_\varPhi \rightarrow S_{\varPhi_v}$ for each $v$ and the diagonal map $\Delta: S_\varPhi \rightarrow S_{\varPhi, \A}$.
For any $\eta = \bigotimes_v \eta_v \in \hat{S}_{\varPhi, \A}$ such that $\eta_v \in \hat{S}_{\varPhi_v, \V_v}$ for all $v$, we may form an irreducible representation $\varSigma_{\eta} = \bigotimes_v \varSigma_{\eta_v}$ of $\SO(\V)(\A)$, where $\varSigma_{\eta_v}$ is the representation of $\SO(\V)(F_v)$ associated to $\eta_v$ by the local Langlands correspondence described in \S \ref{ss:rep-so}.
Then Arthur's multiplicity formula (which has not been established if $\SO(\V)$ is nonsplit but will be assumed in this paper) asserts that
\begin{equation}
\label{eq:mult-so-temp}
 L^2_{\varPhi}(\SO(\V)) \cong \bigoplus_\eta m_{\eta} \varSigma_{\eta},
\end{equation}
where the direct sum runs over continuous characters $\eta$ of $S_{\varPhi, \A}$ such that $\eta_v \in \hat{S}_{\varPhi_v, \V_v}$ for all $v$ and
\[
 m_\eta = 
 \begin{cases}
  1 & \text{if $\Delta^* \eta = \1$;} \\
  0 & \text{otherwise.}
 \end{cases}
\]
In particular, $L^2_{\varPhi}(\SO(\V))$ is multiplicity-free and $m_\disc(\varSigma_\eta) = m_\eta$.

We will use Arthur's multiplicity formula
\begin{itemize}
 \item to globalize a local representation to a global automorphic representation;
 \item to extract local information from a global product formula.
\end{itemize}
To make the argument work, we first globalize a local $L$-parameter to a global $A$-parameter suitably.
For this, we need the following refinement of the globalization of Sakellaridis--Venkatesh \cite[Theorem 16.3.2]{sv}, \cite[Corollary A.8]{ilm}, which will be proved in Appendix \ref{a:global} below.

\begin{prop}
\label{p:global1}
Let $S$ be a nonempty finite set of nonarchimedean places of $\F$ and $v_0$ a nonarchimedean place of $\F$ such that $v_0 \notin S$.
For each $v \in S$, let $\phi_v$ be a $2n$-dimensional irreducible symplectic representation of $L_{\F_v}$.
Then there exists an irreducible symplectic cuspidal automorphic representation $\varPhi$ of $\GL_{2n}(\A)$ such that
\begin{itemize}
 \item $\varPhi_v = \phi_v$ for all $v \in S$;
 \item $\varPhi_v$ is a sum of $1$-dimensional representations for all nonarchimedean $v \notin S \cup \{ v_0 \}$;
 \item $L(\frac{1}{2}, \varPhi) \ne 0$.
\end{itemize}
Here we regard $\varPhi_v$ as a $2n$-dimensional representation of $L_{\F_v}$ via the local Langlands correspondence.
\end{prop}

We also need the following globalization (noting that any irreducible symplectic representation of $L_\R$ is $2$-dimensional), which will be proved in Appendix \ref{a:global} below as well.

\begin{prop}
\label{p:global2}
Assume that $\F$ is totally real.
Let $S$ be a finite set of places of $\F$ containing all archimedean places and $v_0$ a nonarchimedean place of $\F$ such that $v_0 \notin S$.
For each $v \in S$, let $\phi_v$ be a $2$-dimensional irreducible symplectic representation of $L_{\F_v}$.
Then there exists an irreducible symplectic cuspidal automorphic representation $\varPhi$ of $\GL_2(\A)$ such that
\begin{itemize}
 \item $\varPhi_v = \phi_v$ for all $v \in S$;
 \item $\varPhi_v$ is a sum of $1$-dimensional representations for all $v \notin S \cup \{ v_0 \}$;
 \item $L(\frac{1}{2}, \varPhi) \ne 0$.
\end{itemize}
Here we regard $\varPhi_v$ as a $2$-dimensional representation of $L_{\F_v}$ via the local Langlands correspondence.
\end{prop}

Later, we will use the following consequence of these globalizations.

\begin{cor}
\label{c:global}
Let $F$ be a nonarchimedean local field of characteristic zero or $F=\R$.
Let $\phi$ be a $2n$-dimensional good symplectic representation of $L_F$.
\begin{enumerate}
\item \label{c:global1}
Let $\F$ be a totally imaginary number field (resp.~a real quadratic field) if $F$ is nonarchimedean (resp.~if $F=\R$) with two distinct places $v_0, v_1$ of $\F$ such that $\F_{v_0} = \F_{v_1} = F$.
Then there exists a generic elliptic $A$-parameter $\varPhi$ for $\Mp_{2n}$ such that
\begin{itemize}
 \item $\varPhi_{v_0} = \varPhi_{v_1} = \phi$;
 \item the natural maps $S_\varPhi \rightarrow S_{\varPhi_{v_0}}$ and $S_\varPhi \rightarrow S_{\varPhi_{v_1}}$ agree;
 \item $L(\frac{1}{2}, \varPhi) \ne 0$.
\end{itemize}

\item \label{c:global2}
Assume that $\phi$ is reducible.
Let $\F$ be a totally imaginary number field (resp.~$\F=\Q$) if $F$ is nonarchimedean (resp.~if $F=\R$) with a place $v_0$ of $\F$ such that $\F_{v_0} = F$.
Then there exists a generic elliptic $A$-parameter $\varPhi$ for $\Mp_{2n}$ such that
\begin{itemize}
 \item $\varPhi_{v_0} = \phi$;
 \item $\varPhi_v$ is non-good for all $v \ne v_0$;
 \item the natural map $S_{\varPhi} \rightarrow S_{\varPhi_{v_0}}$ is surjective;
 \item the natural map $S_{\varPhi} \rightarrow \prod_{v \ne v_0} S_{\varPhi_v}$ is injective;
 \item $L(\frac{1}{2}, \varPhi) \ne 0$.
\end{itemize}
\end{enumerate}
\end{cor}

\begin{proof}
Put $S_0 = \{ v_0, v_1 \}$ in case \eqref{c:global1}; $S_0 = \{ v_0 \}$ in case \eqref{c:global2}.
We may write $\phi = \bigoplus_i \phi_i$ for some (not necessarily distinct) $2n_i$-dimensional irreducible symplectic representations $\phi_i$ of $L_F$.
For each $i$, choose a set $S_i = \{ v_i, v_i' \}$ consisting of two distinct nonarchimedean places of $\F$ such that
\begin{itemize}
 \item $S_0 \cap S_i = \varnothing$ for all $i$;
 \item $S_i \cap S_j = \varnothing$ for all $i \ne j$.
\end{itemize}
Since we can always find a $2n_i$-dimensional irreducible symplectic representation of $L_{\F_v}$ for all nonarchimedean $v$ (see e.g.~\cite[\S 2]{savin}) and any irreducible representation of $L_\C$ is $1$-dimensional, it follows from Propositions \ref{p:global1} and \ref{p:global2} that there exists an irreducible symplectic cuspidal automorphic representation $\varPhi_i$ of $\GL_{2n_i}(\A)$ such that
\begin{itemize}
 \item $\varPhi_{i,v} = \phi_i$ for all $v \in S_0$;
 \item $\varPhi_{i,v_i}$ is irreducible;
 \item $\varPhi_{i,v}$ is a sum of $1$-dimensional representations for all $v \notin S_0 \cup S_i$;
 \item $L(\frac{1}{2}, \varPhi_i) \ne 0$.
\end{itemize}
In particular, the $\varPhi_i$'s are pairwise distinct.

We show that $\varPhi = \bigoplus_i \varPhi_i$ satisfies the required conditions.
In case \eqref{c:global1}, this is easy.
In case \eqref{c:global2}, since $\phi$ is reducible, $\varPhi_v$ contains a $1$-dimensional irreducible summand (which is non-symplectic) and hence is non-good for all $v \ne v_0$.
Also, since $\varPhi_{i,v_i}$ is irreducible but $\varPhi_{j,v_i}$ is a sum of $1$-dimensional representations for all $j \ne i$, the natural map $S_{\varPhi} \rightarrow \prod_i S_{\varPhi_{v_i}}$ is bijective, so that the natural map $S_\varPhi \rightarrow \prod_{v \ne v_0} S_{\varPhi_v}$ is injective.
The remaining conditions can be easily verified.
\end{proof}

\subsection{The case of good $L$-parameters}
\label{ss:good}

We now prove Propositions \ref{p:moeglin} and \ref{p:key} for any $2n$-dimensional good symplectic representation $\phi$ of $L_F$ when $F$ is nonarchimedean or $F=\R$.
We give a somewhat roundabout argument below, though we could streamline it by using an explicit construction of $A$-packets by M{\oe}glin \cite{m1,m2,m3,m4} and M{\oe}glin--Renard \cite{mr2,mr3}, which were not fully available when this paper was first written.
Especially, when $F$ is nonarchimedean, we could give a local proof by combining M{\oe}glin's results \cite{m1,m2,m3,m4} with a result of Atobe and the first-named author \cite{atobe-gan} which describes local theta lifts explicitly, but we do not discuss it here.

Put $\phi' = \phi \oplus S_{2r-2n}$ with $r > 2n+1$.
For any irreducible representation $\sigma'$ of $\SO_{2r+1}$ and any character $\eta'$ of $\bar{S}_{\phi'}$, we define the multiplicity $m(\sigma', \eta')$ by
\[
 m(\sigma', \eta') = \dim \Hom_{\SO_{2r+1}}(\sigma', \sigma_{\eta'}),
\]
where $\sigma_{\eta'}$ is the representation in the $A$-packet $\Pi_{\phi'}(\SO_{2r+1})$ associated to $\eta'$.

\begin{lem}
\label{l:mult-free}
Let $\sigma'$ be an irreducible representation of $\SO_{2r+1}$.
Then, for any $\eta' \in \hat{\bar{S}}_{\phi'}$, we have
\[
 m(\sigma',\eta') \le 1, 
\]
with equality for at most one $\eta'$.
Namely, $\Pi_{\phi'}(\SO_{2r+1})$ is multiplicity-free.
\end{lem}

\begin{proof}
The assertion was proved in a more general context by M{\oe}glin \cite{m3} when $F$ is nonarchimedean, but is not fully known when $F = \R$, though there has been progress by M{\oe}glin--Renard \cite{mr3}.
Here we give a proof based on theta lifts.

We may assume that $m(\sigma',\eta') > 0$ for some $\eta'$.
Let $\F, v_0, v_1$, and $\varPhi$ be as given in Corollary \ref{c:global}\eqref{c:global1}, so that $\varPhi_{v_0} = \varPhi_{v_1} = \phi$.
Put $\varPhi' = \varPhi \oplus S_{2r-2n}$.
Since $L(\tfrac{1}{2}, \varPhi) \ne 0$, it follows from Lemma \ref{l:epsilon} that the quadratic character $\epsilon_{\varPhi'}$ of $\bar{S}_{\varPhi'}$ is trivial.
We define an abstract irreducible representation $\varSigma' = \bigotimes_v \varSigma'_v$ of $\SO_{2r+1}(\A)$ by setting
\begin{itemize}
 \item $\varSigma'_{v_0} = \varSigma'_{v_1} = \sigma'$;
 \item $\varSigma'_v$ to be the irreducible representation in the $L$-packet $\Pi_{\varphi_{\varPhi'_v}}(\SO_{2r+1})$ associated to the trivial character of $S_{\varphi_{\varPhi'_v}}$ if $v \ne v_0, v_1$.
\end{itemize}
By \cite[Proposition 7.4.1]{a}, $\varSigma'_v$ is an irreducible summand of the representation in the $A$-packet $\Pi_{\varPhi'_v}(\SO_{2r+1})$ associated to the trivial character of $\bar{S}_{\varPhi'_v}$ if $v \ne v_0, v_1$.
Since the pullback of $\eta' \otimes \eta'$ under the natural map 
\[
 \bar{S}_{\varPhi'} \longrightarrow \bar{S}_{\varPhi'_{v_0}} \times \bar{S}_{\varPhi'_{v_1}} = \bar{S}_{\phi'} \times \bar{S}_{\phi'}
\]
is trivial for any $\eta' \in \hat{\bar{S}}_{\phi'}$, we have an embedding
\[
 \Big( \bigoplus_{\eta' \in \hat{\bar{S}}_{\phi'}} ( \sigma_{\eta'} \otimes \sigma_{\eta'} ) \Big)
 \otimes \Big( \bigotimes_{v \ne v_0, v_1} \varSigma'_v \Big) 
 \longhookrightarrow L^2_{\varPhi'}(\SO_{2r+1})
\]
by the multiplicity formula \eqref{eq:mult-so}.
In particular,
\[
 m_\disc(\varSigma') \ge \sum_{\eta' \in \hat{\bar{S}}_{\phi'}} m(\sigma', \eta')^2 > 0.
\]

We now consider theta lifts.
Fix a nontrivial additive character $\varPsi$ of $\F \backslash \A$ such that $\varPsi_{v_0}$ and $\varPsi_{v_1}$ belong to the $(F^\times)^2$-orbit of the fixed nontrivial additive character $\psi$ of $F$.
Let $\varPi = \vartheta_\varPsi^\abs(\varSigma')$ be the abstract theta lift to $\Mp_{2n}(\A)$ relative to $\varPsi$.
Since $\varPhi$ is generic, it follows from Proposition \ref{p:tempered} and Corollary \ref{c:mult-preserve} that $\varPi$ is an irreducible summand of $L^2_{\varPhi,\varPsi}(\Mp_{2n})$ and that 
\[
 m_\cusp(\varPi) = m_\disc(\varPi) = m_\disc(\varSigma') > 0.
\]
For any realization $\mathcal{V} \subset \AA_\cusp(\Mp_{2n})$ of $\varPi$, let $\Theta^\aut_{\V,\varPsi}(\mathcal{V})$ be the global theta lift to $\SO(\V)(\A)$ relative to $\varPsi$, where $\V$ is a $(2n+1)$-dimensional quadratic space over $\F$ with trivial discriminant.
Since $\varPhi$ is generic, we deduce from the tower property and the local theta correspondence for unramified representations that $\Theta^\aut_{\V,\varPsi}(\mathcal{V})$ is cuspidal (but possibly zero).
Moreover, we have
\[
 L^S_\varPsi(\tfrac{1}{2}, \varPi) = L^S(\tfrac{1}{2}, \varPhi) \ne 0,
\]
where $S$ is a sufficiently large finite set of places of $\F$ and $L^S_\varPsi(s,\varPi)$ is the partial $L$-function of $\varPi$ relative to $\varPsi$ and the standard representation of $\Sp_{2n}(\C)$.
Hence, by the Rallis inner product formula \cite{kr,gqt,y} and the local Shimura correspondence \cite{ab1,ab2,gs}, there exists a unique $\V$ such that $\Theta^\aut_{\V,\varPsi}(\mathcal{V})$ is nonzero for any realization $\mathcal{V}$ of $\varPi$.
For this unique $\V$, let $\varSigma = \theta_{\V,\varPsi}^\abs(\varPi)$ be the abstract theta lift to $\SO(\V)(\A)$ relative to $\varPsi$.
Then, by the multiplicity preservation \cite[Proposition 2.6]{g1}, we have
\[
 m_\disc(\varSigma) \ge m_\cusp(\varSigma) = m_\cusp(\varPi).
\]
Also, it follows form the local theta correspondence for unramified representations that $\varSigma$ is an irreducible summand of $L^2_{\varPhi}(\SO(\V))$.
Since $\varPhi$ is generic, the multiplicity formula \eqref{eq:mult-so-temp} implies that 
\[
 m_\disc(\varSigma) = 1.
\]
Thus, combining these (in)equalities, we obtain
\[
 1 \ge \sum_{\eta' \in \hat{\bar{S}}_{\phi'}} m(\sigma', \eta')^2 > 0.
\]
This completes the proof.
\end{proof}

Let $\JH(\Pi_{\phi'}(\SO_{2r+1}))$ be the set of equivalence classes of irreducible representations $\sigma'$ of $\SO_{2r+1}$ such that $m(\sigma', \eta') > 0$ for some $\eta'$.

\begin{lem}
\label{l:packet-inj}
The theta lift induces an injection 
\[
 \vartheta_\psi : \JH(\Pi_{\phi'}(\SO_{2r+1})) \longhookrightarrow \Pi_{\phi,\psi}(\Mp_{2n}).
\]
\end{lem}

\begin{proof}
We retain the notation of the proof of Lemma \ref{l:mult-free}.
In particular, for any $\sigma' \in \JH(\Pi_{\phi'}(\SO_{2r+1}))$, we have an irreducible summand $\varSigma'$ of $L^2_{\varPhi'}(\SO_{2r+1})$ such that
\begin{itemize}
 \item $\varSigma'_{v_0} = \sigma'$;
 \item $\varPi = \vartheta_\varPsi^\abs(\varSigma')$ is an irreducible summand of $L^2_{\varPhi,\varPsi}(\Mp_{2n})$;
 \item $\varSigma = \theta_{\V,\varPsi}^\abs(\varPi)$ is an irreducible summand of $L^2_\varPhi(\SO(\V))$.
\end{itemize}
Since the $L$-parameter of $\varSigma_{v_0} = \theta_{\varPsi_{v_0}}(\varPi_{v_0})$ is $\varPhi_{v_0} = \phi$, we have
\[
 \vartheta_\psi(\sigma') = \vartheta_{\varPsi_{v_0}}(\varSigma'_{v_0}) = \varPi_{v_0} \in \Pi_{\phi,\psi}(\Mp_{2n}).
\]
This and the Howe duality imply the assertion.
\end{proof}

We now show that the map $\vartheta_\psi$ in Lemma \ref{l:packet-inj} is in fact surjective.

\begin{lem}
\label{l:packet-bij}
The theta lift induces a bijection
\[
 \vartheta_\psi : \JH(\Pi_{\phi'}(\SO_{2r+1})) \longleftrightarrow \Pi_{\phi,\psi}(\Mp_{2n}).
\]
\end{lem}

\begin{proof}
Let $\eta' \in \hat{\bar{S}}_{\phi'}$ and put $\eta = \iota^* \eta'$.
Let $\varSigma_\eta$ be the irreducible representation in the $L$-packet $\Pi_\phi(\SO(V))$ associated to $\eta$, where $V$ is the $(2n+1)$-dimensional quadratic space over $F$ with trivial discriminant such that $\eta \in \hat{S}_{\phi,V}$.
Let $\F, v_0, v_1$, and $\varPhi$ be as given in Corollary \ref{c:global}\eqref{c:global1}, so that $\varPhi_{v_0} = \varPhi_{v_1} = \phi$.
Then there exists a unique $(2n+1)$-dimensional quadratic space $\V$ over $\F$ with trivial discriminant such that
\begin{itemize}
 \item $\V_{v_0} = \V_{v_1} = V$;
 \item $\V_v$ is the split space with trivial discriminant for all $v \ne v_0, v_1$.
\end{itemize}
We define an abstract irreducible representation $\varSigma = \bigotimes_v \varSigma_v$ of $\SO(\V)(\A)$ by setting
\begin{itemize}
 \item $\varSigma_{v_0} = \varSigma_{v_1} = \varSigma_\eta$;
 \item $\varSigma_v$ to be the irreducible representation in the $L$-packet $\Pi_{\varPhi_v}(\SO(\V))$ associated to the trivial character of $S_{\varPhi_v}$ if $v \ne v_0, v_1$.
\end{itemize}
Then, by the multiplicity formula \eqref{eq:mult-so-temp}, $\varSigma$ is an irreducible summand of $L^2_{\varPhi}(\SO(\V))$.
In fact, since $\varPhi$ is generic, it follows from the argument in the proof of Proposition \ref{p:tempered} that $\varSigma$ is cuspidal.

Let $\varPi = \theta_\varPsi^\abs(\varSigma)$ be the abstract theta lift to $\Mp_{2n}(\A)$ (relative to $\varPsi$ as in the proof of Lemma \ref{l:mult-free}), so that 
\[
 \varPi_{v_0} = \theta_{\varPsi_{v_0}}(\varSigma_{v_0}) = \theta_\psi(\varSigma_\eta) = \pi_\eta \in \Pi_{\phi,\psi}(\Mp_{2n}).
\]
Since $L(\frac{1}{2}, \varPhi) \ne 0$, we can deduce from the argument in the proof of Lemma \ref{l:mult-free} that the global theta lift $\Theta_\varPsi^\aut(\varSigma)$ to $\Mp_{2n}(\A)$ is nonzero and cuspidal, so that $\Theta_\varPsi^\aut(\varSigma) \cong \varPi$ is an irreducible summand of $L^2_{\varPhi,\varPsi}(\Mp_{2n})$.
Hence, by Corollary \ref{c:mult-preserve}, the abstract theta lift $\varSigma' = \vartheta_\varPsi^\abs(\varPi)$ to $\SO_{2r+1}(\A)$ is an irreducible summand of $L^2_{\varPhi'}(\SO_{2r+1})$, where $\varPhi' = \varPhi \oplus S_{2r-2n}$.
Since $\varPhi'_{v_0} = \phi'$, this implies that
\[
 \varSigma'_{v_0} \in \JH(\Pi_{\phi'}(\SO_{2r+1})).
\]
On the other hand, we have
\[
 \vartheta_{\varPsi_{v_0}}(\varSigma'_{v_0}) = \varPi_{v_0} = \pi_\eta.
\]
Hence the map in Lemma \ref{l:packet-inj} is surjective.
\end{proof}

Finally, we show that Propositions \ref{p:moeglin} and \ref{p:key} hold for $\phi$.
We consider the irreducible case and the reducible case separately.

\begin{lem}
Assume that $\phi$ is irreducible.
Then Propositions \ref{p:moeglin} and \ref{p:key} hold for $\phi$.
\end{lem}

\begin{proof}
Let $\sigma'_\1$ be the representation in the $A$-packet $\Pi_{\phi'}(\SO_{2r+1})$ associated to the trivial character of $\bar{S}_{\phi'}$.
Then, by \cite[Proposition 7.4.1]{a}, $\sigma'_\1$ contains the irreducible representation $\sigma'$ in the $L$-packet $\Pi_{\varphi_{\phi'}}(\SO_{2r+1})$ associated to the trivial character of $S_{\varphi_{\phi'}}$.
Since 
\[
 \# \JH(\Pi_{\phi'}(\SO_{2r+1})) = 2
\]
by the irreducibility of $\phi$ and Lemma \ref{l:packet-bij}, we may write $\JH(\Pi_{\phi'}(\SO_{2r+1})) = \{ \sigma', \sigma'' \}$ for some irreducible representation $\sigma''$ of $\SO_{2r+1}$.

We show that $\sigma'_\1$ is irreducible.
Suppose on the contrary that $\sigma'_\1$ is reducible.
By Lemma \ref{l:mult-free}, we have $\sigma'_\1 = \sigma' \oplus \sigma''$.
Let $\F, v_0, v_1$, and $\varPhi$ be as given in Corollary \ref{c:global}\eqref{c:global1}, so that $\varPhi_{v_0} = \varPhi_{v_1} = \phi$.
Put $\varPhi' = \varPhi \oplus S_{2r-2n}$.
Since $L(\tfrac{1}{2}, \varPhi) \ne 0$, it follows from Lemma \ref{l:epsilon} that the quadratic character $\epsilon_{\varPhi'}$ of $\bar{S}_{\varPhi'}$ is trivial.
We define an abstract irreducible representation $\varSigma' = \bigotimes_v \varSigma'_v$ of $\SO_{2r+1}(\A)$ by setting
\begin{itemize}
 \item $\varSigma'_{v_0} = \sigma'$;
 \item $\varSigma'_{v_1} = \sigma''$;
 \item $\varSigma'_v$ to be the irreducible representation in the $L$-packet $\Pi_{\varphi_{\varPhi'_v}}(\SO_{2r+1})$ associated to the trivial character of $S_{\varphi_{\varPhi'_v}}$ if $v \ne v_0, v_1$.
\end{itemize}
Since $\sigma'_\1 = \sigma' \oplus \sigma''$, it follows from the multiplicity formula \eqref{eq:mult-so} that $\varSigma'$ is an irreducible summand of $L^2_{\varPhi'}(\SO_{2r+1})$.
Let $\varPi = \vartheta_\varPsi^\abs(\varSigma')$ be the abstract theta lift to $\Mp_{2n}(\A)$ (relative to $\varPsi$ as in the proof of Lemma \ref{l:mult-free}).
Then, as shown in the proof of Lemma \ref{l:mult-free}, there exists a unique $(2n+1)$-dimensional quadratic space $\V$ over $\F$ with trivial discriminant such that the abstract theta lift $\theta^\abs_{\V,\varPsi}(\varPi)$ to $\SO(\V)(\A)$ is nonzero.
On the other hand, if $v \ne v_0, v_1$ (so that $v$ is not real), then by Lemma \ref{l:theta-triv-char}, $\varPi_v = \vartheta_{\varPsi_v}(\varSigma'_v)$ is the irreducible representation in the $L$-packet $\Pi_{\varPhi_v,\varPsi_v}(\Mp_{2n})$ associated to the trivial character of $S_{\varPhi_v}$, so that
\[
 \varepsilon(\V_v) = 1.
\]
Also, by Lemma \ref{l:packet-bij}, we have
\[
 \{ \varPi_{v_0}, \varPi_{v_1} \} = \{ \vartheta_{\varPsi_{v_0}}(\varSigma'_{v_0}), \vartheta_{\varPsi_{v_1}}(\varSigma'_{v_1}) \} = \{ \vartheta_\psi(\sigma'), \vartheta_\psi(\sigma'') \} = \Pi_{\phi,\psi}(\Mp_{2n}), 
\]
so that 
\[
 \varepsilon(\V_{v_0}) \cdot \varepsilon(\V_{v_1}) = -1.
\]
This contradicts the fact that $\prod_v \varepsilon(\V_v) = 1$.
Hence $\sigma'_\1$ is irreducible and $\sigma'_\1 = \sigma'$.

By Lemma \ref{l:mult-free} and \ref{l:packet-bij}, it remains to show that $\vartheta_\psi(\sigma')$ is associated to the trivial character of $S_\phi$.
If $F$ is nonarchimedean, then this follows from Lemma \ref{l:theta-triv-char}.
Suppose that $F = \R$ (so that $n=1$) and that $\vartheta_\psi(\sigma')$ is associated to the nontrivial character of $S_\phi$.
From the above argument with the following modifications:
\begin{itemize}
 \item $\F = \Q$;
 \item $\varPhi$ is a generic elliptic $A$-parameter for $\Mp_2$ such that $\varPhi_\infty = \phi$ and $L(\frac{1}{2}, \varPhi) \ne 0$ (see Proposition \ref{p:global2});
 \item $\varSigma'_\infty = \sigma'$;
 \item $\varSigma'_v$ is associated to the trivial character of $S_{\varphi_{\varPhi'_v}}$ if $v$ is nonarchimedean,
\end{itemize}
we can deduce that there exists a $3$-dimensional quadratic space $\V$ over $\F$ with trivial discriminant such that $\varepsilon(\V_v) = 1$ for all nonarchimedean $v$ but $\varepsilon(\V_\infty) = -1$.
This is a contradiction and completes the proof.
\end{proof}

\begin{lem}
Assume that $\phi$ is reducible.
Then Propositions \ref{p:moeglin} and \ref{p:key} hold for $\phi$.
\end{lem}

\begin{proof}
Let $\eta' \in \hat{\bar{S}}_{\phi'}$ and put $\eta = \iota^* \eta'$.
Let $\F, v_0$, and $\varPhi$ be as given in Corollary \ref{c:global}\eqref{c:global2}, so that $\varPhi_{v_0} = \phi$.
Since the natural map $(\prod_{v \ne v_0} S_{\varPhi_v}) \hat{\mathstrut} \rightarrow \hat{S}_{\varPhi}$ is surjective, there exists a continuous character $\bigotimes_v \xi_v$ of $S_{\varPhi, \A}$ such that
\begin{align}
 \xi_{v_0} & = \eta, \label{eq:xi_v0} \\
 \Big( \bigotimes_v \xi_v \Big) \circ \Delta & = \1. \label{eq:prod1}
\end{align}
Then there exists a unique $(2n+1)$-dimensional quadratic space $\V$ over $\F$ with trivial discriminant such that $\xi_v \in \hat{S}_{\varPhi_v, \V_v}$ for all $v$.
We define an abstract irreducible representation $\varSigma = \bigotimes_v \varSigma_v$ of $\SO(\V)(\A)$ by setting
\[
 \varSigma_v = \varSigma_{\xi_v} \in \Pi_{\varPhi_v}(\SO(\V))
\]
for all $v$. 
Then, by the multiplicity formula \eqref{eq:mult-so-temp}, $\varSigma$ is an irreducible summand of $L^2_{\varPhi}(\SO(\V))$.

Let $\varPi = \theta_\varPsi^\abs(\varSigma)$ be the abstract theta lift to $\Mp_{2n}(\A)$ (relative to a fixed nontrivial additive character $\varPsi$ of $\F \backslash \A$ such that $\varPsi_{v_0}$ belongs to the $(F^\times)^2$-orbit of $\psi$), so that 
\[
 \varPi_v = \theta_{\varPsi_v}(\varSigma_v) = \theta_{\varPsi_v}(\varSigma_{\xi_v}) = \pi_{\xi_v} \in \Pi_{\varPhi_v,\varPsi_v}(\Mp_{2n})
\]
for all $v$.
Then, as shown in the proof of Lemma \ref{l:packet-bij}, the abstract theta lift $\varSigma' = \vartheta_\varPsi^\abs(\varPi)$ to $\SO_{2r+1}(\A)$ is an irreducible summand of $L^2_{\varPhi'}(\SO_{2r+1})$, where $\varPhi' = \varPhi \oplus S_{2r-2n}$.
This implies that for any $v$, $\varSigma'_v$ is an irreducible summand of $\sigma_{\xi'_v} \in \Pi_{\varPhi_v'}(\SO_{2r+1})$ for some $\xi'_v \in \hat{\bar{S}}_{\varPhi'_v}$.
By Lemma \ref{l:mult-free}, such $\xi_v'$ is unique.
In fact, if $v \ne v_0$ (so that $v$ is not real), then since $\varPhi_v$ is non-good, we can apply Propositions \ref{p:moeglin} and \ref{p:key} to obtain $\varSigma'_v = \sigma_{\xi'_v}$ with 
\begin{equation}
\label{eq:eta_v}
 \iota_v^* \xi'_v = \xi_v,
\end{equation}
where $\iota_v : S_{\varPhi_v} \rightarrow \bar{S}_{\varPhi'_v}$ is the natural isomorphism.
On the other hand, since $L(\tfrac{1}{2}, \varPhi) \ne 0$, it follows from Lemma \ref{l:epsilon} that the quadratic character $\epsilon_{\varPhi'}$ of $\bar{S}_{\varPhi'}$ is trivial.
Hence, by the multiplicity formula \eqref{eq:mult-so}, we must have
\begin{equation}
\label{eq:prod2}
 \Big( \bigotimes_v \xi'_v \Big) \circ \Delta = \1.
\end{equation}
Thus, since the natural map $\hat{S}_{\varPhi_{v_0}} \rightarrow \hat{S}_{\varPhi}$ is injective, we deduce from \eqref{eq:prod1}, \eqref{eq:eta_v}, \eqref{eq:prod2} that
\[
 \iota_{v_0}^* \xi'_{v_0} = \xi_{v_0}.
\]
Hence we have $\iota_{v_0}^* \xi'_{v_0} = \eta$ by \eqref{eq:xi_v0}, so that
\begin{equation}
\label{eq:eta_v0}
 \xi'_{v_0} = \eta'
\end{equation}
by the definition of $\eta$.
In particular, $\sigma_{\eta'} \in \Pi_{\phi'}(\SO_{2r+1})$ is nonzero for any $\eta' \in \hat{\bar{S}}_{\phi'}$.
By Lemmas \ref{l:mult-free} and \ref{l:packet-bij}, this implies that $\sigma_{\eta'}$ is irreducible for any $\eta'$.
Moreover, by \eqref{eq:eta_v0}, we have
\[
 \vartheta_\psi(\sigma_{\eta'}) = \vartheta_{\varPsi_{v_0}}(\varSigma'_{v_0}) = \varPi_{v_0} = \pi_\eta.
\]
This completes the proof.
\end{proof}

This completes the proof of Propositions \ref{p:moeglin} and \ref{p:key} and hence of Theorem \ref{t:main2}.

\appendix

\section{\textbf{Globalizations}}
\label{a:global}

In this appendix, we prove Propositions \ref{p:global1} and \ref{p:global2}.

\subsection{Proof of Proposition \ref{p:global1}}

Let $F$ be a number field and $\A$ the ad\`ele ring of $F$.
Let $S_\infty$ be the set of archimedean places of $F$.
Let $S$ be a nonempty finite set of nonarchimedean places of $F$ and $v_0$ a nonarchimedean place of $F$ such that $v_0 \notin S$.
Then Proposition \ref{p:global1} is an immediate consequence of the following:

\begin{prop}
\label{a:global-gl1}
For each $v \in S \cup \{ v_0 \}$, let $\tau_v$ be an irreducible square-integrable representation of $\GL_{2n}(F_v)$ such that $L(s, \tau_v, \wedge^2)$ has a pole at $s=0$.
Assume that $\tau_{v_0}$ is supercuspidal.
Then there exists an irreducible cuspidal automorphic representation $\TT$ of $\GL_{2n}(\A)$ such that 
\begin{itemize}
 \item $\TT_v = \tau_v$ for all $v \in S \cup \{ v_0 \}$;
 \item $\TT_v$ is a principal series representation for all $v \notin S_\infty \cup S \cup \{ v_0 \}$;
 \item $L(s, \TT, \wedge^2)$ has a pole at $s = 1$;
 \item $L(\frac{1}{2}, \TT) \ne 0$.
\end{itemize}
\end{prop}

To prove this proposition, we globalize a generic representation of $\Mp_{2n}(F_v)$ to a globally generic automorphic representation of $\Mp_{2n}(\A)$.
More precisely, let $N$ be the standard maximal unipotent subgroup of $\Sp_{2n}$.
We regard $N(\A)$ as a subgroup of $\Mp_{2n}(\A)$ via the canonical splitting.
Fix a nontrivial additive character $\psi$ of $F \backslash \A$.
As in \cite[\S 12]{ggp}, $\psi$ gives rise to a generic character of $N(F) \backslash N(\A)$, which we denote again by $\psi$.

\begin{prop}
\label{a:global-mp}
For each $v \in S \cup \{ v_0 \}$, let $\pi_v$ be an irreducible genuine $\psi_v$-generic square-integrable representation of $\Mp_{2n}(F_v)$.
Assume that $\pi_{v_0}$ is supercuspidal, and that if $n=1$, then $\pi_{v_0}$ is not the odd Weil representation relative to $\psi_{v_0}$.
Then there exists an irreducible genuine globally $\psi$-generic cuspidal automorphic representation $\varPi$ of $\Mp_{2n}(\A)$ such that
\begin{itemize}
 \item $\varPi_v = \pi_v$ for all $v \in S \cup \{ v_0 \}$;
 \item $\varPi_v$ is a principal series representation for all $v \notin S_\infty \cup S \cup \{ v_0 \}$.
\end{itemize}
\end{prop}

Proposition \ref{a:global-gl1} can be easily deduced from Proposition \ref{a:global-mp} and \cite[Proposition 4.3]{ilm}.
We include the proof for the convenience of the reader.

\begin{proof}[Proof of Proposition \ref{a:global-gl1}]
For each $v \in S \cup \{ v_0 \}$, let $\pi_v$ be the descent of $\tau_v$ to $\Mp_{2n}(F_v)$ relative to $\psi_v$ (see \cite{grs99}, \cite{grs02}, \cite[Theorem 3.1]{ilm}).
Then $\pi_v$ satisfies the conditions in Proposition \ref{a:global-mp}.
Let $\varPi$ be as given in Proposition \ref{a:global-mp}.
Let $\varSigma = \Theta_\psi^\aut(\varPi)$ be the global theta lift to $\SO_{2n+1}(\A)$ relative to $\psi$.
Since $\theta_{\psi_{v_0}}(\pi_{v_0})$ is supercuspidal by \cite[Theorem 2.2]{js03}, $\varSigma$ is cuspidal.
By \cite{furusawa}, $\varSigma$ is nonzero and globally generic.
Hence $\varSigma$ is irreducible, and by the Rallis inner product formula \cite{kr,gqt,y}, we have
\[
 L_\psi(\tfrac{1}{2}, \varPi) \ne 0,
\]
where $L_\psi(s, \varPi)$ is the $L$-function of $\varPi$ relative to $\psi$ and the standard representation of $\Sp_{2n}(\C)$.
Moreover, it follows from the local Shimura correspondence \cite{gs} that
\begin{itemize}
 \item $\varSigma_v$ is square-integrable for all $v \in S \cup \{ v_0 \}$;
 \item $\varSigma_v$ is a principal series representation for all $v \notin S_\infty \cup S \cup \{ v_0 \}$.
\end{itemize}

We now take $\TT$ to be the functorial lift of $\varSigma$ to $\GL_{2n}(\A)$.
Then $L(s, \TT, \wedge^2)$ has a pole at $s=1$, and we have
\[
 L(\tfrac{1}{2}, \TT) = L(\tfrac{1}{2}, \varSigma) = L_\psi(\tfrac{1}{2}, \varPi) \ne 0.
\]
Also, since $\TT_v$ is the functorial lift of $\varSigma_v$, it is a principal series representation for all $v \notin S_\infty \cup S \cup \{ v_0 \}$.
Finally, by \cite[Proposition 4.3]{ilm}, we have
\[
 \TT_v = \tau_v
\]
for all $v \in S \cup \{ v_0 \}$.
In particular, $\TT_{v_0}$ is supercuspidal and hence $\TT$ is cuspidal.
\end{proof}

It remains to prove Proposition \ref{a:global-mp}, which is a refinement of the globalization of Sakellaridis--Venkatesh \cite[Theorem 16.3.2]{sv}, \cite[Corollary A.8]{ilm}.
We need to modify their argument to control the localization $\varPi_v$ at nonarchimedean $v$ outside $S \cup \{ v_0 \}$.

\begin{proof}[Proof of Proposition \ref{a:global-mp}]
We first introduce some notation.
By abuse of notation, we write
\[
 G_v = \Mp_{2n}(F_v), \qquad
 G_S = \Mp_{2n}(F_S), \qquad 
 G(\A) = \Mp_{2n}(\A), 
\]
where $F_S = \prod_{v \in S} F_v$.
If $F_v$ is nonarchimedean of odd residual characteristic, we regard $K_v = \Sp_{2n}(\mathcal{O}_v)$ as a subgroup of $G_v$ via the standard splitting, where $\mathcal{O}_v$ is the integer ring of $F_v$.
Also, we regard $G(F) = \Sp_{2n}(F)$ as a subgroup of $G(\A)$ via the canonical splitting.
Put 
\[
 N_S = N(F_S), \qquad
 \psi_S = \bigotimes_{v \in S} \psi_v, \qquad
 \pi_S = \bigotimes_{v \in S} \pi_v.
\]
Let $C_c^\infty(G_v)$ be the space of genuine smooth functions on $G_v$ with compact support.
Let $C_c^\infty(N_v \backslash G_v, \psi_v)$ be the space of genuine smooth functions $f$ on $G_v$ such that
\begin{itemize}
 \item $\supp f$ is compact modulo $N_v$;
 \item $f(xg) = \psi_v(x) f(g)$ for all $x \in N_v$ and $g \in G_v$.
\end{itemize}
Then we have a map $\mathcal{P}_v :C_c^\infty(G_v) \rightarrow C_c^\infty(N_v \backslash G_v, \psi_v)$ defined by
\[
 (\mathcal{P}_v \tilde{f})(g) = \int_{N_v} \tilde{f}(xg) \overline{\psi_v(x)} \, dx.
\]
Also, $C_c^\infty(N_v \backslash G_v, \psi_v)$ is equipped with a natural hermitian inner product $\langle \cdot, \cdot \rangle$.
Let $C_c^\infty(N_S \backslash G_S, \psi_S)$ be defined similarly.
For any automorphic form $\varphi$ on $G(\A)$, we define its Whittaker--Fourier coefficient $\mathcal{W}_\varphi$ by 
\[
 \mathcal{W}_\varphi(g) = \int_{N(F) \backslash N(\A)} \varphi(x g) \overline{\psi(x)} \, dx.
\]
Let $L^2(G)$ be the genuine part of $L^2(G(F) \backslash G(\A))$, which is equipped with the Petersson inner product $\langle \cdot, \cdot \rangle$.
Let $L^2_\cusp(G)$ be the closure in $L^2(G)$ of the subspace of genuine cusp forms on $G(\A)$.
We define $L^2_{\cusp, \psi\text{-}\mathrm{gen}}(G)$ as the orthogonal complement in $L^2_\cusp(G)$ of the closure of the subspace of genuine cusp forms $\varphi$ on $G(\A)$ such that $\mathcal{W}_\varphi = 0$.
Then, by the uniqueness of Whittaker models, $L^2_{\cusp, \psi\text{-}\mathrm{gen}}(G)$ is multiplicity-free.
Fix a finite set $S_0$ of nonarchimedean places of $F$ such that 
\begin{itemize}
 \item $S \cap S_0 = \varnothing$;
 \item $v_0 \notin S_0$;
 \item $F_v$ is of odd residual characteristic and $\psi_v$ is of order zero if $v \notin S_\infty \cup S \cup \{ v_0 \} \cup S_0$.
\end{itemize}
Let $\{ \varPi_i \}$ be the set of irreducible summands of $L^2_{\cusp, \psi\text{-}\mathrm{gen}}(G)$ such that
\begin{itemize}
 \item $\varPi_{i,{v_0}} = \pi_{v_0}$;
 \item $\varPi_{i,v}$ is a principal series representation for all $v \in S_0$;
 \item $\varPi_{i,v}$ is unramified (i.e.~$\varPi_{i,v}$ has a nonzero $K_v$-fixed vector) for all $v \notin S_\infty \cup S \cup \{ v_0 \} \cup S_0$.
\end{itemize}

As in \cite[\S 16.4]{sv}, we show that $\pi_S$ is weakly contained in the Hilbert space direct sum $\bigoplus_i \varPi_i$ (regarded as a representation of $G_S$).
Since $\pi_S$ belongs to the support of the Plancherel measure for $L^2(N_S \backslash G_S, \psi_S)$, it suffices to prove the following: for any compact subset $\Omega \subset G_S$ and any $f_S \in C_c^\infty(N_S \backslash G_S, \psi_S)$ with $\langle f_S, f_S \rangle = 1$, there exist $c_i \in \R$ with $c_i \ge 0$ and $\varphi_i \in \varPi_i$ with $\langle \varphi_i, \varphi_i \rangle = 1$ such that $\sum_i c_i = 1$ and 
\[
 \langle R_y f_S, f_S \rangle = \sum_i c_i \langle R_y \varphi_i, \varphi_i \rangle
\]
for all $y \in \Omega$.
Here $R_y$ is the right translation by $y$.
For each $v \notin S_\infty \cup S$, we choose $f_v = \mathcal{P}_v \tilde{f}_v \in C_c^\infty(N_v \backslash G_v, \psi_v)$, where $\tilde{f}_v \in C_c^\infty(G_v)$ is given as follows:
\begin{itemize}
 \item if $v = v_0$, then $\tilde{f}_{v_0}$ is a matrix coefficient of $\pi_{v_0}$;
 \item if $v \in S_0$, then $\tilde{f}_v$ belongs to the wave packet associated to a Bernstein component consisting only of irreducible principal series representations (e.g.~one containing a sufficiently ramified principal series representation);
 \item if $v \notin S_\infty \cup S \cup \{ v_0 \} \cup S_0$, then
 \[
  \tilde{f}_v(g) = 
  \begin{cases}
   \epsilon & \text{if $g \in \epsilon \cdot K_v$ with $\epsilon \in \{ \pm 1 \}$;} \\
   0 & \text{otherwise.}
  \end{cases}
 \]
\end{itemize}
By \cite[Lemma 4.4]{psp}, we may further assume that $\langle f_v, f_v \rangle = 1$.
For each $v \in S_\infty$, we choose $f_v \in C_c^\infty(N_v \backslash G_v, \psi_v)$ later.
Put
\[
 f = f_S \otimes \Big( \bigotimes_{v \notin S} f_v \Big).
\]
We define a Poincar\'e series $P_f$ on $G(\A)$ by 
\[
 P_f(g) = \sum_{\gamma \in N(F) \backslash G(F)} f(\gamma g),
\]
where the sum converges absolutely.
Then we have
\begin{align*}
 \langle P_f, P_f \rangle 
 & \le \int_{N(F) \backslash G(\A)} |f(g) P_f(g)| \, dg \\
 & \le \int_{N(\A) \backslash G(\A)} |f(g)| \int_{N(F) \backslash N(\A)} |P_f(xg)| \, dx \, dg < \infty,
\end{align*}
so that $P_f$ is square-integrable over $G(F) \backslash G(\A)$.
Since $\pi_{v_0}$ is supercuspidal, $P_f$ is cuspidal by the condition on $f_{v_0}$.
For any genuine cusp form $\varphi$ on $G(\A)$ such that $\mathcal{W}_\varphi = 0$, we have
\begin{align*}
 \langle P_f, \varphi \rangle
 & = \int_{N(F) \backslash G(\A)} f(g) \overline{\varphi(g)} \, dg \\
 & = \int_{N(\A) \backslash G(\A)} f(g) \overline{\mathcal{W}_\varphi(g)} \, dg = 0. 
\end{align*}
Hence $P_f \in L^2_{\cusp, \psi\text{-}\mathrm{gen}}(G)$.
Moreover, by the conditions on $f_v$, we have $P_f \in \bigoplus_i \varPi_i$.
Thus, as in \cite[(16.8)]{sv}, it remains to prove the following:
there exists $f_v \in C_c^\infty(N_v \backslash G_v, \psi_v)$ for $v \in S_\infty$ such that
\[
 \langle R_y f_S, f_S \rangle = \langle R_y P_f, P_f \rangle 
\]
for all $y \in \Omega$.
This was proved in \cite[\S 16.4]{sv}, but we include the proof for the convenience of the reader.
We choose $f_v \in C_c^\infty(N_v \backslash G_v, \psi_v)$ for $v \in S_\infty$ such that $\langle f_v, f_v \rangle = 1$ and $\supp f_v$ is sufficiently small, so that there exists a compact subset $U \subset G(\A)$ satisfying the following conditions:
\begin{itemize}
 \item $U^{-1} = U$;
 \item $\Omega \subset U$;
 \item $\supp R_y f \subset N(\A) \cdot U$ for all $y \in \Omega$;
 \item $N(\A) \subset N(F) \cdot U$;
 \item $G(F) \cap U_4 \cdot N(\A) = N(F)$, where $U_4 = U \cdot U \cdot U \cdot U$.
\end{itemize}
For any $y \in \Omega$, we have
\begin{align*}
 \langle R_y P_f, P_f \rangle 
 & = \int_{N(F) \backslash G(\A)} f(gy) \overline{P_f(g)} \, dg \\
 & = \int_{N(F) \backslash G(\A)} \sum_{\gamma \in N(F) \backslash G(F)}
 f(gy) \overline{f(\gamma g)} \, dg.
\end{align*}
Let $\gamma \in G(F)$.
If $f(gy) \overline{f(\gamma g)} \ne 0$ for some $g \in G(\A)$, then $gy \in N(\A) \cdot U$ and $\gamma g \in N(\A) \cdot U \subset N(F) \cdot U \cdot U$.
Replacing $\gamma$ by an element in $N(F) \cdot \gamma$ if necessary, we may assume that $\gamma g \in U \cdot U$.
Then we have $\gamma = \gamma g \cdot y \cdot (g y)^{-1} \in U_4 \cdot N(\A)$, so that $\gamma \in N(F)$.
Hence 
\[
 \langle R_y P_f, P_f \rangle
 = \int_{N(F) \backslash G(\A)} f(g y) \overline{f(g)} \, dg
 = \langle R_y f, f \rangle
 = \langle R_y f_S, f_S \rangle.
\]

Thus, we have shown that $\pi_S$ is weakly contained in $\bigoplus_i \varPi_i$.
As shown in the proof of Proposition \ref{a:global-gl1}, the global theta lift $\Theta_\psi^\aut(\varPi_i)$ to $\SO_{2n+1}(\A)$ is cuspidal.
Hence, by \cite[Proposition A.7]{ilm}, we have
\[
 \{ \varPi_{i,S} \} \subset \Irr_{\mathrm{unit}, \psi_S\text{-}\mathrm{gen}, \le c} G_S
\]
for some $0 \le c < \frac{1}{2}$, where $\varPi_{i,S} = \bigotimes_{v \in S} \varPi_{i,v}$ and $\Irr_{\mathrm{unit}, \psi_S\text{-}\mathrm{gen}, \le c} G_S$ is the set of equivalence classes of irreducible genuine $\psi_S$-generic unitary representations of $G_S$ whose exponents are bounded by $c$ (see \cite[p.~1328]{ilm}).
On the other hand, it follows from the analog of \cite[Lemma A.2]{ilm} for $\Mp_{2n}$ that $\pi_S$ is isolated in $\Irr_{\mathrm{unit}, \psi_S\text{-}\mathrm{gen}, \le c} G_S$ with respect to the Fell topology.
This forces $\pi_S = \varPi_{i,S}$ for some $i$ and completes the proof, noting that $\varPi_{i,v}$ is an irreducible principal series representation as in Remark \ref{r:unitarity} for all $v \notin S_\infty \cup S \cup \{ v_0 \} \cup S_0$.
\end{proof}

\subsection{Proof of Proposition \ref{p:global2}}

Let $F$ be a totally real number field and $\A$ the ad\`ele ring of $F$.
Let $S_\infty$ be the set of archimedean places of $F$.
Let $S$ be a finite set of places of $F$ containing $S_\infty$ and $v_0$ a nonarchimedean place of $F$ such that $v_0 \notin S$.
Then Proposition \ref{p:global2} is an immediate consequence of the following:

\begin{prop}
\label{a:global-gl2}
For each $v \in S$, let $\tau_v$ be an irreducible square-integrable representation of $\GL_2(F_v)$ with trivial central character.
Then there exists an irreducible cuspidal automorphic representation $\TT$ of $\GL_2(\A)$ with trivial central character such that
\begin{itemize}
\item $\TT_v = \tau_v$ for all $v \in S$;
\item $\TT_v$ is a principal series representation for all $v \notin S \cup \{ v_0 \}$;
\item $L(\frac{1}{2}, \TT) \ne 0$.
\end{itemize}
\end{prop}

This proposition can be deduced from Waldspurger's result \cite[Th\'eor\`eme 4]{w2} on the nonvanishing of central $L$-values (see also \cite{fh}, in which Friedberg--Hoffstein gave a completely different proof).
Nevertheless, we give a more direct proof based on torus periods, so that the proof of the main result of this paper will be independent of Waldspurger's nonvanishing result.

We write 
\[
 S = S_\infty \cup S_1 \cup S_2,
\]
where $S_1$ (resp.~$S_2$) is the set of nonarchimedean places $v$ in $S$ such that $\tau_v$ is not supercuspidal (resp.~$\tau_v$ is supercuspidal).
Put 
\[
 S_0 = 
\begin{cases}
 S_\infty \cup S_1 & \text{if $\# (S_\infty \cup S_1)$ is even}; \\
 S_\infty \cup S_1 \cup \{ v_0 \} & \text{if $\# (S_\infty \cup S_1)$ is odd}.
\end{cases}
\]
Let $D$ be the quaternion division algebra over $F$ which is ramified precisely at the places in $S_0$.
For each $v \in S$, let $\pi_v$ be the Jacquet--Langlands transfer of $\tau_v$ to $D_v^\times$.
More explicitly, we have:
\begin{itemize}
\item if $v \in S_\infty$, then $\pi_v$ is odd-dimensional;
\item if $v \in S_1$, then $\pi_v = \chi_v \circ \nu_v$ for some quadratic character $\chi_v$ of $F_v^\times$, where $\nu_v$ is the reduced norm on $D_v$;
\item if $v \in S_2$, then $\pi_v$ is supercuspidal.
\end{itemize}
Choose a totally imaginary quadratic extension $E$ of $F$ such that
\begin{itemize}
\item if $v \in S_1$ and $\chi_v$ is unramified, then $E_v$ is the unramified quadratic extension of $F_v$;
\item if $v \in S_1$ and $\chi_v$ is ramified, then $E_v$ is the quadratic extension of $F_v$ associated to $\chi_v$ by class field theory;
\item if $v \in S_2$, then $E_v = F_v \times F_v$;
\item if $v = v_0$, then $E_v \ne F_v \times F_v$.
\end{itemize}
Then $E$ can be embedded into $D$.
Fix an embedding $E \hookrightarrow D$.

We consider algebraic groups
\[
 G = D^\times/F^\times, \qquad
 H = E^\times/F^\times
\]
defined over $F$, and regard $H$ as a subgroup of $G$.
Note that $G$ is anisotropic over $F$ and $G_v$ is compact for all $v \in S_0$.
Then, for any $v \in S$, $\pi_v$ (regarded as a representation of $G_v$) is $H$-distinguished, i.e.~
\begin{equation}
\label{eq:local-H-dist}
 \Hom_{H_v}(\pi_v, \C) \ne 0. 
\end{equation}
We globalize these distinguished representations to a globally distinguished automorphic representation.

\begin{prop}
\label{a:global-gp}
There exists an irreducible automorphic representation $\varPi$ of $G(\A)$ such that
\begin{itemize}
\item $\varPi_v = \pi_v$ for all $v \in S$;
\item $\varPi_v$ is a principal series representation for all $v \notin S \cup \{ v_0 \}$;
\item $\varPi$ is not $1$-dimensional;
\item $\varPi$ is globally $H$-distinguished, i.e.~
\[
 \int_{H(F) \backslash H(\A)} \varphi(h) \, dh \ne 0
\]
for some $\varphi \in \varPi$.
\end{itemize}
\end{prop}

Proposition \ref{a:global-gl2} can be easily deduced from Proposition \ref{a:global-gp} and Waldspurger's result \cite{wal} on torus periods.
Indeed, let $\varPi$ be as given in Proposition \ref{a:global-gp} and $\TT$ the Jacquet--Langlands transfer of $\varPi$ to $\GL_2(\A)$.
Since $\varPi$ is not $1$-dimensional, $\TT$ is cuspidal.
Moreover, since $\varPi$ is globally $H$-distinguished, it follows from Waldspurger's result \cite[Th\'eor\`eme 2]{wal} (see also \cite{jacquet-w}, in which Jacquet gave a new proof based on relative trace formulas) that 
\[
 L(\tfrac{1}{2}, \TT) \cdot L(\tfrac{1}{2}, \TT \times \chi_E) \ne 0,
\]
where $\chi_E$ is the quadratic automorphic character of $\A^\times$ associated to $E/F$ by class field theory.
Hence $\TT$ satisfies the required conditions.

Proposition \ref{a:global-gp} follows from the globalization of Prasad--Schulze-Pillot \cite[Theorem 4.1]{psp}.
We include the proof for the convenience of the reader.

\begin{proof}[Proof of Proposition \ref{a:global-gp}]
We first introduce some notation.
If $v \notin S_0$, we denote by $K_v = \PGL_2(\mathcal{O}_v)$ the standard maximal compact subgroup of $G_v = \PGL_2(F_v)$, where $\mathcal{O}_v$ is the integer ring of $F_v$.
Let $C_c^\infty(G_v)$ and $C_c^\infty(H_v \backslash G_v)$ be the spaces of smooth functions on $G_v$ and $H_v \backslash G_v$ with compact support, respectively. 
Then we have a map $\mathcal{P}_v :C_c^\infty(G_v) \rightarrow C_c^\infty(H_v \backslash G_v)$ defined by
\[
 (\mathcal{P}_v \tilde{f})(g) = \int_{H_v} \tilde{f}(hg) \, dh.
\]
Put $L^2(G) = L^2(G(F) \backslash G(\A))$, which is equipped with the Petersson inner product $\langle \cdot, \cdot \rangle$.
We define $L^2_0(G)$ as the orthogonal complement in $L^2(G)$ of the closure of the subspace spanned by all automorphic characters of $G(\A)$.

For each $v \ne v_0$, we choose $f_v = \mathcal{P}_v \tilde{f}_v \in C_c^\infty(H_v \backslash G_v)$, where $\tilde{f}_v \in C_c^\infty(G_v)$ is given as follows:
\begin{itemize}
\item if $v \in S$, then $\tilde{f}_v$ is a matrix coefficient of $\pi_v$;
\item if $v \notin S \cup \{ v_0 \}$, then $\tilde{f}_v$ is the characteristic function of $K_v$.
\end{itemize}
By \eqref{eq:local-H-dist} and \cite[Lemma 4.4]{psp} (which continues to hold for $v \in S_\infty$), we may further assume that $f_v(1) \ne 0$.
For $v = v_0$, we choose $f_{v_0} \in C_c^\infty(H_{v_0} \backslash G_{v_0})$ later.
Put $f = \bigotimes_v f_v$.
We define a Poincar\'e series $P_f$ on $G(\A)$ by 
\[
 P_f(g) = \sum_{\gamma \in H(F) \backslash G(F)} f(\gamma g),
\]
where the sum converges absolutely.
Obviously, $P_f \in L^2(G)$.
Also, for any automorphic form $\varphi$ on $G(\A)$, we have
\begin{equation}
\label{eq:poincare}
 \langle P_f, \varphi \rangle = \int_{H(F) \backslash G(\A)} f(g) \overline{\varphi(g)} \, dg = \int_{H(\A) \backslash G(\A)} f(g) \int_{H(F) \backslash H(\A)} \overline{\varphi(hg)} \, dh \, dg.
\end{equation}
We choose $f_{v_0} \in C_c^\infty(H_{v_0} \backslash G_{v_0})$ such that $\supp f_{v_0}$ is sufficiently small, so that
\begin{itemize}
\item $\supp f \cap G(F) \subset H(F)$;
\item $\supp f_{v_0} \subset H_{v_0} \cdot \ker \mu$ for all characters $\mu$ of $G_{v_0}$, 
\end{itemize}
and such that $f_{v_0}(1) \ne 0$ and
\[
 \int_{H_{v_0} \backslash G_{v_0}} f_{v_0}(g) \, dg = 0.
\]
Then we have
\[
 P_f(1) = f(1) \ne 0.
\]
Let $\mu$ be an automorphic character of $G(\A)$.
By \eqref{eq:poincare}, we have $\langle P_f, \mu \rangle = 0$ unless $\mu$ is trivial on $H(\A)$, in which case we have
\begin{align*}
 \langle P_f, \mu \rangle & = \operatorname{vol}(H(F) \backslash H(\A)) \cdot \int_{H(\A) \backslash G(\A)} f(g) \overline{\mu(g)} \, dg \\
 & = \operatorname{vol}(H(F) \backslash H(\A)) \cdot \prod_v \int_{H_v \backslash G_v} f_v(g_v) \overline{\mu_v(g_v)} \, dg_v = 0
\end{align*}
as well since 
\[
 \int_{H_{v_0} \backslash G_{v_0}} f_{v_0}(g) \overline{\mu_{v_0}(g)} \, dg = \int_{H_{v_0} \backslash G_{v_0}} f_{v_0}(g) \, dg = 0.
\]
Hence we have $P_f \in L^2_0(G)$.

We now take $\varPi$ to be any irreducible automorphic representation of $G(\A)$ such that $\langle P_f, \varphi \rangle \ne 0$ for some $\varphi \in \varPi$.
In particular, $\varPi$ is not $1$-dimensional.
Moreover, by \eqref{eq:poincare}, we have
\[
 \int_{H(F) \backslash H(\A)} \varphi(hg) \, dh \ne 0 
\]
for some $g \in G(\A)$, so that $\varPi$ is globally $H$-distinguished.
Finally, by the conditions on $f_v$, we have: 
\begin{itemize}
\item $\varPi_v = \pi_v$ for all $v \in S$;
\item $\varPi_v$ is unramified (i.e.~$\varPi_v$ has a nonzero $K_v$-fixed vector) for all $v \notin S \cup \{ v_0 \}$.
\end{itemize}
Hence $\varPi$ satisfies the required conditions.
\end{proof}


\begin{thebibliography}{99}

\bibitem{adams}
J.~Adams,
\emph{Lifting of characters on orthogonal and metaplectic groups.}
Duke Math. J. \textbf{92} (1998), no.~1, 129--178.

\bibitem{ab1}
J.~Adams and D.~Barbasch,
\emph{Reductive dual pair correspondence for complex groups.}
J. Funct. Anal. \textbf{132} (1995), no.~1, 1--42.

\bibitem{ab2}
J.~Adams and D.~Barbasch,
\emph{Genuine representations of the metaplectic group.}
Compos. Math. \textbf{113} (1998), no.~1, 23--66.

\bibitem{amr}
N.~Arancibia, C.~M{\oe}glin, and D.~Renard,
\emph{Paquets d'Arthur des groupes classiques et unitaires.}
arXiv:1507.01432

\bibitem{a2011}
J.~Arthur,
\emph{The embedded eigenvalue problem for classical groups.}
On certain $L$-functions, 19--31, Clay Math. Proc. \textbf{13}, Amer. Math. Soc., Providence, RI, 2011.

\bibitem{a}
J.~Arthur,
\emph{The endoscopic classification of representations: orthogonal and symplectic groups.}
American Mathematical Society Colloquium Publications \textbf{61},
American Mathematical Society, Providence, RI, 2013.

\bibitem{atobe}
H.~Atobe,
\emph{The local theta correspondence and the local Gan--Gross--Prasad conjecture for the symplectic-metaplectic case.}
Math. Ann. \textbf{371} (2018), no.~1-2, 225--295.

\bibitem{atobe-gan}
H.~Atobe and W.~T.~Gan,
\emph{Local theta correspondence of tempered representations and Langlands parameters.}
Invent. Math. \textbf{210} (2017), no.~2, 341--415.

\bibitem{bj}
D.~Ban and C.~Jantzen,
\emph{The Langlands quotient theorem for finite central extensions of $p$-adic groups.}
Glas. Mat. Ser. III \textbf{48(68)} (2013), no.~2, 313--334.

\bibitem{bergeron-clozel}
N.~Bergeron and L.~Clozel,
\emph{Quelques cons\'equences des travaux d'Arthur pour le spectre et la topologie des vari\'et\'es hyperboliques.} 
Invent. Math. \textbf{192} (2013), no.~3, 505--532.

\bibitem{bf}
D.~Bump and S.~Friedberg,
\emph{The exterior square automorphic $L$-functions on $\mathrm{GL}(n)$.}
Festschrift in honor of I.~I.~Piatetski-Shapiro on the occasion of his sixtieth birthday, Part II, 47--65, Israel Math. Conf. Proc. \textbf{3}, Weizmann, Jerusalem, 1990.

\bibitem{bg} 
D.~Bump and D.~Ginzburg,
\emph{Symmetric square $L$-functions on $\mathrm{GL}(r)$.}
Ann. of Math. (2) \textbf{136} (1992), no.~1, 137--205.

\bibitem{flicker} 
Y.~Z.~Flicker, 
\emph{Automorphic forms on covering groups of $\mathrm{GL}(2)$.}
Invent. Math. \textbf{57} (1980), no.~2, 119--182.

\bibitem{fh} 
S.~Friedberg and J.~Hoffstein,
\emph{Nonvanishing theorems for automorphic $L$-functions on $\mathrm{GL}(2)$.}
Ann. of Math. (2) \textbf{142} (1995), no.~2, 385--423.

\bibitem{furusawa}
M.~Furusawa,
\emph{On the theta lift from $\mathrm{SO}_{2n+1}$ to $\widetilde{\mathrm{Sp}}_n$.}
J. Reine Angew. Math. \textbf{466} (1995), 87--110.

\bibitem{g1}
W.~T.~Gan,
\emph{The Saito--Kurokawa space of $\mathrm{PGSp}_4$ and its transfer to inner forms.}
Eisenstein series and applications, 87--123, Progr. Math. \textbf{258}, Birkh{\"a}user Boston, Boston, MA, 2008.

\bibitem{g2}
W.~T.~Gan,
\emph{A Langlands program for covering groups?}
Proceedings of the Sixth International Congress of Chinese Mathematicians, Vol.~I, 57--78, Adv. Lect. Math. \textbf{36}, Int. Press, Somerville, MA, 2017.

\bibitem{g3}
W.~T.~Gan,
\emph{The Shimura correspondence \`a la Waldspurger.}
\url{http://www.math.nus.edu.sg/~matgwt/postech.pdf}

\bibitem{ggp}
W.~T.~Gan, B.~H.~Gross, and D.~Prasad,
\emph{Symplectic local root numbers, central critical $L$-values, and restriction problems in the representation theory of classical groups.}
Sur les conjectures de Gross et Prasad. I,
Ast\'erisque \textbf{346} (2012), 1--109.

\bibitem{gi1}
W.~T.~Gan and A.~Ichino,
\emph{Formal degrees and local theta correspondence.}
Invent. Math. \textbf{195} (2014), no.~3, 509--672.

\bibitem{gi2}
W.~T.~Gan and A.~Ichino,
\emph{The Gross--Prasad conjecture and local theta correspondence.}
Invent. Math. \textbf{206} (2016), no.~3, 705--799.

\bibitem{gi-mp-real}
W.~T.~Gan and A.~Ichino,
\emph{On the irreducibility of some induced representations of real reductive Lie groups.}
Tunisian J. Math. \textbf{1} (2019), no.~1, 73--107.

\bibitem{gqt}
W.~T.~Gan, Y.~Qiu, and S.~Takeda,
\emph{The regularized Siegel--Weil formula (the second term identity) and the Rallis inner product formula.}
Invent. Math. \textbf{198} (2014), no.~3, 739--831.

\bibitem{gs}
W.~T.~Gan and G.~Savin,
\emph{Representations of metaplectic groups I: epsilon dichotomy and local Langlands correspondence.}
Compos. Math. \textbf{148} (2012), no.~6, 1655--1694.

\bibitem{gt1}
W.~T.~Gan and S.~Takeda,
\emph{On the Howe duality conjecture in classical theta correspondence.}
Advances in the theory of automorphic forms and their $L$-functions, 105--117, Contemp. Math. \textbf{664}, Amer. Math. Soc., Providence, RI, 2016.

\bibitem{gt2}
W.~T.~Gan and S.~Takeda,
\emph{A proof of the Howe duality conjecture.}
J. Amer. Math. Soc. \textbf{29} (2016), no.~2, 473--493.

\bibitem{grs99}
D.~Ginzburg, S.~Rallis, and D.~Soudry,
\emph{On a correspondence between cuspidal representations of $\mathrm{GL}_{2n}$ and $\widetilde{\mathrm{Sp}}_{2n}$.}
J. Amer. Math. Soc. \textbf{12} (1999), no.~3, 849--907.

\bibitem{grs02}
D.~Ginzburg, S.~Rallis, and D.~Soudry,
\emph{Endoscopic representations of $\widetilde{\mathrm{Sp}}_{2n}$.}
J. Inst. Math. Jussieu \textbf{1} (2002), no.~1, 77--123.

\bibitem{gp}
B.~H.~Gross and D.~Prasad,
\emph{On the decomposition of a representation of $\mathrm{SO}_n$ when restricted to $\mathrm{SO}_{n-1}$.}
Canad. J. Math. \textbf{44} (1992), no.~5, 974--1002.

\bibitem{ht}
M.~Harris and R.~Taylor,
\emph{The geometry and cohomology of some simple Shimura varieties.}
Annals of Mathematics Studies \textbf{151},
Princeton University Press, Princeton, NJ, 2001.

\bibitem{henniart}
G.~Henniart,
\emph{Une preuve simple des conjectures de Langlands pour $\mathrm{GL}(n)$ sur un corps $p$-adique.}
Invent. Math. \textbf{139} (2000), no.~2, 439--455.

\bibitem{howe0}
R.~Howe,
\emph{$\theta$-series and invariant theory.}
Automorphic forms, representations and $L$-functions, Part 1, 275--285,
Proc. Sympos. Pure Math. \textbf{33}, Amer. Math. Soc., Providence, RI, 1979.

\bibitem{howe1}
R.~Howe,
\emph{Automorphic forms of low rank.}
Noncommutative harmonic analysis and Lie groups, 211--248,
Lecture Notes in Math. \textbf{880}, Springer, Berlin-New York, 1981.

\bibitem{howe2}
R.~Howe,
\emph{On a notion of rank for unitary representations of the classical groups.}
Harmonic analysis and group representations, 223--331, Liguori, Naples, 1982.

\bibitem{howe3}
R.~Howe,
\emph{Transcending classical invariant theory.}
J. Amer. Math. Soc. \textbf{2} (1989), no.~3, 535--552.

\bibitem{ilm}
A.~Ichino, E.~Lapid, and Z.~Mao,
\emph{On the formal degrees of square-integrable representations of odd special orthogonal and metaplectic groups.}
Duke Math. J. \textbf{166} (2017), no.~7, 1301--1348.

\bibitem{jacquet-w}
H.~Jacquet, 
\emph{Sur un r\'esultat de Waldspurger.}
Ann. Sci. \'Ecole Norm. Sup. (4) \textbf{19} (1986), no.~2, 185--229.

\bibitem{jacquet}
H.~Jacquet,
\emph{On the nonvanishing of some $L$-functions.}
Proc. Indian Acad. Sci. Math. Sci. \textbf{97} (1987), no.~1-3, 117--155.

\bibitem{js}
H.~Jacquet and J.~A.~Shalika,
\emph{On Euler products and the classification of automorphic forms. II.}
Amer. J. Math. \textbf{103} (1981), no.~4, 777--815.

\bibitem{js-aa}
H.~Jacquet and J.~A.~Shalika,
\emph{Exterior square $L$-functions.}
Automorphic forms, Shimura varieties, and $L$-functions, Vol.~II, 143--226, Perspect. Math. \textbf{11}, Academic Press, Boston, MA, 1990.

\bibitem{js03}
D.~Jiang and D.~Soudry,
\emph{The local converse theorem for $\mathrm{SO}(2n+1)$ and applications.}
Ann. of Math. (2) \textbf{157} (2003), no.~3, 743--806.

\bibitem{kudla}
S.~S.~Kudla,
\emph{On the local theta-correspondence.}
Invent. Math. \textbf{83} (1986), no.~2, 229--255.

\bibitem{kr}
S.~S.~Kudla and S.~Rallis,
\emph{A regularized Siegel--Weil formula: the first term identity.}
Ann. of Math. (2) \textbf{140} (1994), no.~1, 1--80.

\bibitem{langlands1}
R.~P.~Langlands,
\emph{On the notion of an automorphic representation.}
Automorphic forms, representations and $L$-functions, Part 1, 203--207,
Proc. Sympos. Pure Math. \textbf{33}, Amer. Math. Soc., Providence, RI, 1979.

\bibitem{langlands2}
R.~P.~Langlands,
\emph{On the classification of irreducible representations of real algebraic groups.}
Representation theory and harmonic analysis on semisimple Lie groups, 101--170, Math. Surveys Monogr. \textbf{31}, Amer. Math. Soc., Providence, RI, 1989.

\bibitem{lmt}
E.~Lapid, G.~Mui\'c, and M.~Tadi\'c,
\emph{On the generic unitary dual of quasisplit classical groups.}
Int. Math. Res. Not. \textbf{2004}, no.~26, 1335--1354.

\bibitem{li1}
J.-S.~Li,
\emph{Singular unitary representations of classical groups.}
Invent. Math. \textbf{97} (1989), no.~2, 237--255.

\bibitem{li2}
J.-S.~Li,
\emph{On the classification of irreducible low rank unitary representations of classical groups.}
Compos. Math. \textbf{71} (1989), no.~1, 29--48.

\bibitem{li3}
J.-S.~Li,
\emph{Automorphic forms with degenerate Fourier coefficients.}
Amer. J. Math. \textbf{119} (1997), no.~3, 523--578.

\bibitem{wwli1}
W.-W.~Li,
\emph{Transfert d'int\'egrales orbitales pour le groupe m\'etaplectique.}
Compos. Math. \textbf{147} (2011), no.~2, 524--590.

\bibitem{wwli2}
W.-W.~Li,
\emph{La formule des traces stable pour le groupe m\'etaplectique: les termes elliptiques.}
Invent. Math. \textbf{202} (2015), no.~2, 743--838.

\bibitem{wwli3}
W.-W.~Li,
\emph{Spectral transfer for metaplectic groups. I. Local character relations.}
To appear in J. Inst. Math. Jussieu, DOI:10.1017/S1474748016000384.

\bibitem{luo}
C.~Luo,
\emph{Endoscopic character identities for metaplectic groups.}
arXiv:1801.10302

\bibitem{m1}
C.~M{\oe}glin,
\emph{Sur certains paquets d'Arthur et involution d'Aubert--Schneider--Stuhler g\'en\'eralis\'ee.}
Represent. Theory \textbf{10} (2006), 86--129.

\bibitem{m2}
C.~M{\oe}glin,
\emph{Paquets d'Arthur discrets pour un groupe classique $p$-adique.}
Automorphic forms and $L$-functions II. Local aspects, 179--257, Contemp. Math. \textbf{489}, Amer. Math. Soc., Providence, RI, 2009.

\bibitem{m3}
C.~M{\oe}glin,
\emph{Multiplicit\'e $1$ dans les paquets d'Arthur aux places $p$-adiques.}
On certain $L$-functions, 333--374, Clay Math. Proc. \textbf{13}, Amer. Math. Soc., Providence, RI, 2011.

\bibitem{m4}
C.~M{\oe}glin,
\emph{Image des op\'erateurs d'entrelacements normalis\'es et p\^oles des s\'eries d'Eisenstein.}
Adv. Math. \textbf{228} (2011), no.~2, 1068--1134.

\bibitem{m5}
C.~M{\oe}glin,
\emph{Paquets d'Arthur sp\'eciaux unipotents aux places archim\'ediennes et correspondance de Howe.}
Representation theory, number theory, and invariant theory, 469--502, Progr. Math. \textbf{323}, Birkh\"auser/Springer, Cham, 2017.

\bibitem{mr}
C.~M{\oe}glin and D.~Renard,
\emph{Paquets d'Arthur des groupes classiques complexes.}
Around Langlands correspondences, 203--256, Contemp. Math. \textbf{691}, Amer. Math. Soc., Providence, RI, 2017.

\bibitem{mr2}
C.~M{\oe}glin and D.~Renard,
\emph{Sur les paquets d'Arthur des groupes classiques r\'eels.}
arXiv:1703.07226

\bibitem{mr3}
C.~M{\oe}glin and D.~Renard,
\emph{Sur les paquets d'Arthur aux places r\'eelles, translation.}
arXiv:1704.05096

\bibitem{mr4}
C.~M{\oe}glin and D.~Renard,
\emph{Sur les paquets d'Arthur des groupes classiques et unitaires non quasi-d\'eploy\'es.}
arXiv:1803.07662

\bibitem{mvw}
C.~M{\oe}glin, M.-F.~Vign\'eras, and J.-L.~Waldspurger,
\emph{Correspondances de Howe sur un corps $p$-adique.}
Lecture Notes in Math. \textbf{1291}, Springer-Verlag, Berlin, 1987.

\bibitem{mw}
C.~M{\oe}glin and J.-L.~Waldspurger,
\emph{La conjecture locale de Gross--Prasad pour les groupes sp\'eciaux orthogonaux: le cas g\'en\'eral.}
Sur les conjectures de Gross et Prasad. II,
Ast\'erisque \textbf{347} (2012), 167--216.

\bibitem{mw1}
C.~M{\oe}glin and J.-L.~Waldspurger,
\emph{Stabilisation de la formule des traces tordue, Volume 1.}
Progr. Math. \textbf{316}, Birkh\"auser/Springer, Cham, 2016.

\bibitem{mw2}
C.~M{\oe}glin and J.-L.~Waldspurger,
\emph{Stabilisation de la formule des traces tordue, Volume 2.}
Progr. Math. \textbf{317}, Birkh\"auser/Springer, Cham, 2016.

\bibitem{ngo}
B.~C.~Ng\^o,
\emph{Le lemme fondamental pour les alg\`ebres de Lie.}
Publ. Math. Inst. Hautes \'Etudes Sci. \textbf{111}, (2010), 1--169.

\bibitem{niwa}
S.~Niwa,
\emph{Modular forms of half integral weight and the integral of certain theta-functions.}
Nagoya Math. J. \textbf{56} (1975), 147--161.

\bibitem{ps}
I.~I.~Piatetski-Shapiro,
\emph{Work of Waldspurger.}
Lie group representations, II, 280--302, Lecture Notes in Math. \textbf{1041}, Springer, Berlin, 1984.

\bibitem{psp}
D.~Prasad and R.~Schulze-Pillot,
\emph{Generalised form of a conjecture of Jacquet and a local consequence.}
J. Reine Angew. Math. \textbf{616} (2008), 219--236.

\bibitem{rangarao}
R.~Ranga Rao,
\emph{On some explicit formulas in the theory of Weil representation.}
Pacific J. Math. \textbf{157} (1993), no.~2, 335--371.

\bibitem{renard}
D.~Renard,
\emph{Endoscopy for $\mathrm{Mp}(2n,\mathbb{R})$.}
Amer. J. Math. \textbf{121} (1999), no.~6, 1215--1243.

\bibitem{sv}
Y.~Sakellaridis and A.~Venkatesh,
\emph{Periods and harmonic analysis on spherical varieties.}
Ast\'erisque \textbf{396} (2017).

\bibitem{savin}
G.~Savin,
\emph{Lifting of generic depth zero representations of classical groups.}
J. Algebra \textbf{319} (2008), no.~8, 3244--3258.

\bibitem{scharlau}
W.~Scharlau,
\emph{Quadratic and Hermitian forms.}
Grundlehren der Mathematischen Wissenschaften \textbf{270}, Springer-Verlag, Berlin, 1985.

\bibitem{scholze}
P.~Scholze,
\emph{The local Langlands correspondence for $\mathrm{GL}_n$ over $p$-adic fields.}
Invent. Math. \textbf{192} (2013), no.~3, 663--715.

\bibitem{shahidi1}
F.~Shahidi,
\emph{On certain $L$-functions.}
Amer. J. Math. \textbf{103} (1981), no.~2, 297--355.

\bibitem{shahidi2}
F.~Shahidi,
\emph{A proof of Langlands' conjecture on Plancherel measures; complementary series for $p$-adic groups.}
Ann. of Math. (2) \textbf{132} (1990), no.~2, 273--330.

\bibitem{shelstad1}
D.~Shelstad,
\emph{$L$-indistinguishability for real groups.}
Math. Ann. \textbf{259} (1982), no.~3, 385--430.

\bibitem{shelstad2}
D.~Shelstad,
\emph{Tempered endoscopy for real groups. III. Inversion of transfer and $L$-packet structure.}
Represent. Theory \textbf{12} (2008), 369--402.

\bibitem{shimura} 
G.~Shimura,
\emph{On modular forms of half integral weight.} 
Ann. of Math. (2) \textbf{97} (1973), 440--481. 

\bibitem{shintani}
T.~Shintani,
\emph{On construction of holomorphic cusp forms of half integral weight.}
Nagoya Math. J. \textbf{58} (1975), 83--126.

\bibitem{spehv}
B.~Speh and D.~A.~Vogan, Jr.,
\emph{Reducibility of generalized principal series representations.}
Acta Math. \textbf{145} (1980), no.~3-4, 227--299.

\bibitem{sun}
B.~Sun,
\emph{Dual pairs and contragredients of irreducible representations.}
Pacific J. Math. \textbf{249} (2011), no.~2, 485--494.

\bibitem{sz}
B.~Sun and C.-B.~Zhu,
\emph{Conservation relations for local theta correspondence.}
J. Amer. Math. Soc. \textbf{28} (2015), no.~4, 939--983.

\bibitem{tadic}
M.~Tadi\'c,
\emph{Classification of unitary representations in irreducible representations of general linear group (non-Archimedean case).}
Ann. Sci. \'Ecole Norm. Sup. (4) \textbf{19} (1986), no.~3, 335--382.

\bibitem{vogan81}
D.~A.~Vogan,~Jr.,
\emph{Representations of real reductive Lie groups.}
Progr. Math. \textbf{15}, Birkh\"auser, Boston, MA, 1981.

\bibitem{w1}
J.-L.~Waldspurger,
\emph{Correspondance de Shimura.}
J. Math. Pures Appl. (9) \textbf{59} (1980), no.~1, 1--132.

\bibitem{wal}
J.-L.~Waldspurger,
\emph{Sur les valeurs de certaines fonctions $L$ automorphes en leur centre de sym\'etrie.}
Compos. Math. \textbf{54} (1985), no.~2, 173--242.

\bibitem{w3}
J.-L.~Waldspurger,
\emph{D\'emonstration d'une conjecture de dualit\'e de Howe dans le cas $p$-adique, $p \neq 2$.}
Festschrift in honor of I.~I.~Piatetski-Shapiro on the occasion of his sixtieth birthday, Part I, 267--324, Israel Math. Conf. Proc. \textbf{2}, Weizmann, Jerusalem, 1990.

\bibitem{w2}
J.-L.~Waldspurger,
\emph{Correspondances de Shimura et quaternions.}
Forum Math. \textbf{3} (1991), no.~3, 219--307.

\bibitem{xu}
B.~Xu,
\emph{On M{\oe}glin's parametrization of Arthur packets for $p$-adic quasisplit $Sp(N)$ and $SO(N)$.}
Canad. J. Math. \textbf{69} (2017), no.~4, 890--960.

\bibitem{y}
S.~Yamana,
\emph{$L$-functions and theta correspondence for classical groups.}
Invent. Math. \textbf{196} (2014), no.~3, 651--732.

\end{thebibliography}
\end{document}